\documentclass[11pt]{article}
\usepackage{times}

\usepackage{bm}
\usepackage{graphicx}

\usepackage{amsfonts}
\usepackage{amsmath}

\newcommand{\bn}{{\mbox{\boldmath$\nu$}}}

\begin{document}
\title{Asymptotic distribution of a consistent cross-spectrum estimator
based on uniformly spaced samples of a non-bandlimited process}

\author{Radhendushka Srivastava and~Debasis Sengupta\\[2ex] Applied Statistics Unit\\ Indian
Statistical Institute\\ 203 B.T. Road\\ Kolkata 700108, India\\
(e-mail: radhe\_r@isical.ac.in; sdebasis@isical.ac.in).}

\maketitle \setcounter{page}{0} \thispagestyle{empty}
\newpage
\begin{center}
{\LARGE\bf Asymptotic distribution of a consistent cross-spectrum
estimator based on uniformly spaced samples of a non-bandlimited
process }
\end{center}

\section*{Abstract}
It is well known that if the power spectral density
of a continuous time stationary stochastic process does not have a
compact support, data sampled from that process at any uniform
sampling rate leads to biased and inconsistent spectrum estimators.
In a recent paper, the authors showed that the smoothed periodogram
estimator can be consistent, if the sampling interval is allowed to
shrink to zero at a suitable rate as the sample size goes to
infinity. In this paper, this `shrinking asymptotics' approach is
used to obtain the limiting distribution of the smoothed periodogram
estimator of spectra and cross-spectra. It is shown that, under
suitable conditions, the scaling that ensures weak convergence of
the estimator to a limiting normal random vector can range from
cube-root of the sample size to square-root of the sample size,
depending on the strength of the assumption made. The results are
used to construct asymptotic confidence intervals for spectra and
cross spectra. It is shown through a Monte-Carlo simulation study
that these intervals have appropriate empirical coverage
probabilities at moderate sample sizes.

{\it Keywords:} Power spectral density, spectrum estimation,
smoothed periodogram, shrinking asymptotics, asymptotic confidence
interval.

\section{Introduction}\label{s1}
Estimation of power spectral density (spectrum) of a continuous
time, mean square continuous, stationary stochastic process is a
classical problem. Generally the estimation is based on finitely
many samples of the process. It is well known that if the spectrum
is compactly supported (bandlimited), then it can be estimated from
uniformly spaced samples, provided the sampling is done at the
Nyquist rate or faster (Kay, 1999). For sampled non-bandlimited
processes, or bandlimited processes sampled at sub-Nyquist rate, the
problem of aliasing leads to biased estimation. For this reason, it
is sometimes argued that a non-bandlimited spectrum can never be
estimated consistently from uniformly spaced samples at any fixed
sampling rate (Shapiro and Silverman, 1960; Masry, 1978).

Consequently, some researchers have turned to non-uniform sampling
schemes such as stochastic sampling and periodic non-uniform
sampling. Masry (1978) proved the consistency of some spectrum
estimators based on stochastic sampling schemes, under appropriate
conditions that allow the underlying spectrum to be non-bandlimited.
For bandlimited processes, it has been shown that periodic
non-uniform sampling at sub-Nyquist average rate can lead to
consistent spectrum estimation (Marvasti, 2001).

It is important to note that the argument of inconsistency of
spectrum estimators computed from uniformly spaced samples is based
on the assumption that the sampling rate remains fixed even as the
sample size goes to infinity. However, when one has the resources to
increase the sample size indefinitely, one would like to use some of
those resources to sample faster, rather than being constrained by a
fixed sampling rate. In fact, it has been shown (Srivastava and
Sengupta, 2010) that if this constraint is removed, and the sampling
rate is allowed to increase suitably as the sample size goes to
infinity, then the smoothed periodogram can be a consistent
estimator of a non-bandlimited spectral density.

It should be noted that uniform sampling is generally far easier to
implement than non-uniform sampling. For this reason, the fact of
consistency of a spectrum estimator computed from uniformly spaced
samples of a non-bandlimited process is noteworthy. This fact also
gives rise to further questions about the estimator, such as its
convergence in distribution, and construction of asymptotic
confidence intervals for the spectrum, based on the estimator. These
questions have so far not been addressed through asymptotic
calculations that allow the sampling rate to go to infinity. This is
what we propose to do in this paper.

The asymptotic approach chosen here (referred to as `shrinking
asymptotics' by Fuentes, 2002) was also adopted by other authors
(e.g., Constantine and Hall, 1994; Hall et al., 1994; Lahiri, 1999),
although the asymptotic distribution of the smoothed periodogram has
not been studied previously. This approach is different from the
`fixed-domain asymptotics' or `infill asymptotics' approach (Chen et
al., 2000; Stein, 1995; Zhang and Zimmerman, 2005; Lim and Stein,
2008) which, in the present case, would have required that the
time-span of the original continuous-time data (before sampling)
remains fixed as the sampling rate goes to infinity.

Let ${\bf X}=\{{\bf X}(t),~-\infty<t<\infty\}$ be a vector-valued
mean square continuous stationary stochastic process, having zero
mean. We denote the components of the process $\bf{X}$ by
$X_a=\{X_a(t),~-\infty<t<\infty \}$ for $a\in\{1,2,\ldots,r\}$, and
the variance-covariance matrix of the process $\bf{X}$ at lag $\tau$
by
\[ \bf{C}(\tau)=\left( \begin{array}{cccc}
C_{11}(\tau) & C_{12}(\tau) &\ldots &C_{1r}(\tau) \\
C_{21}(\tau) & C_{22}(\tau) &\ldots &C_{2r}(\tau) \\
\vdots           &\vdots    &  &\vdots             \\
C_{r1}(\tau) & C_{r2}(\tau)   &\ldots &C_{rr}(\tau) \end{array}
\right),\] where
$$C_{a_1a_2}(\tau)=E\left[X_{a_1}(t+\tau)X_{a_2}(t)\right]~\mbox{for}~a_1,a_2\in\{1,2,\ldots,r\}.$$
The spectral and cross-spectral density matrix of the process
$\bf{X}$ is denoted by
\[ \bf{\Phi}(\cdot)=\left( \begin{array}{cccc}
\phi_{11}(\cdot) & \phi_{12}(\cdot) &\ldots &\phi_{1r}(\cdot) \\
\phi_{21}(\cdot) & \phi_{22}(\cdot) &\ldots &\phi_{2r}(\cdot) \\
\vdots           &\vdots            &\vdots &\vdots           \\
\phi_{r1}(\cdot) & \phi_{r2}(\cdot) &\ldots &\phi_{rr}(\cdot)
\end{array} \right)\] %
where
$$\phi_{a_1a_2}(\lambda)=\frac1{2\pi}\int_{-\infty}^{\infty}C_{a_1a_2}(t)e^{-it\lambda}dt
\quad\mbox{for}~-\infty<\lambda<\infty, \quad
a_1,a_2\in\{1,2,\ldots,r\}.$$

In this paper, we consider the following estimator of
$\phi_{a_1a_2}(\lambda)$ for $a_1~,a_2\in\{1,2,\ldots,r\}$:
\begin{equation}\label{est}
\widehat{\phi}_{a_1a_2}(\lambda)= \frac{1}{2\pi
n\rho_n}\sum_{t_1=1}^{n}\sum_{t_2=1}^{n}K(b_n(t_1-t_2))X_{a_1}
\left(\frac{t_1}{\rho_n}\right)X_{a_2}\left(\frac{t_2}{\rho_n}\right)e^{-\frac{i(t_1-t_2)\lambda}{\rho_n}}
1_{[-\pi\rho_n,\pi\rho_n]}(\lambda),
\end{equation}
where $K(\cdot)$ is a covariance averaging kernel, $b_n$ is the
kernel bandwidth, $\rho_n$ is the sampling rate and $1_A(\lambda)$
is the indicator of the event $\lambda\in A$.

In Section 2, we establish the consistency of the spectrum estimator
(\ref{est}) for non-bandlimited processes, which is a generalization
of a result of Srivastava and Sengupta (2010) to the case of
multivariate time series. It paves the way for our main result on
the asymptotic distribution of the estimator, given in Section~3.
Section~4 contains some discussion on optimal rates of convergence.
In Section~5, we investigate the question as to how large the sample
size should be in order for the applicability of the asymptotic
distribution and the resulting pointwise confidence intervals. We
look for answers through a Monte Carlo simulation study and report
the findings. All the proofs are given in the appendix.

\section{Consistency}
In order to establish the consistency of the estimator
$\widehat{\phi}_{a_1a_2}(\cdot)$ given in (\ref{est}), we make a few
assumptions on the process ${\bf X}$, the kernel $K(\cdot)$ and the
sequences $b_n$ and $\rho_n$.

\bigskip
{\sc Assumption 1.} The function $g_{a_1a_2}(\cdot)$, defined
over the real line as $g_{a_1a_2}(t)=\sup_{|s|\ge
|t|}|C_{a_1a_2}(s)|$ is integrable
for all $a_1,a_2\in\{1,2,\ldots,r\}$.%
\bigskip

{\sc Assumption 2.} The covariance averaging kernel function
$K(\cdot)$ is continuous, even, 
square integrable and bounded by a non-negative, even and integrable
function having a unique maximum at 0. Further, $K(0)=1$.%

\bigskip
{\sc Assumption 3.} The kernel window width is such that
$nb_n\rightarrow\infty$ as $n\rightarrow\infty$.%

\bigskip
{\sc Assumption 4.} The sampling rate is such that
$\rho_n\rightarrow\infty$ and $\rho_n b_n\rightarrow 0$ as
$n\rightarrow \infty$.%
\bigskip

Note that Assumption~4 implies that $b_n\rightarrow 0$ as
$n\rightarrow\infty$.

\bigskip{\sc Theorem~1.} {\em Under Assumptions 1--4, the bias of the
estimator $\widehat{\phi}_{a_1a_2}(\cdot)$ tends to zero uniformly
over any closed and finite interval.}
\bigskip

In order to establish convergence of the variance-covariance matrix,
we need a further assumption involving cumulants. Recall that the
$r$-th order joint cumulant of the random variable
$(Y_1,\ldots,Y_r)$ is given by
\begin{equation}
cum(Y_1,\ldots,Y_r)=\sum_{\bn}
(-1)^{p-1}(p-1)!\left(E\prod_{j\in~\nu_{1}}Y_j\right)\times\cdots
\times\left(E\prod_{j\in~\nu_{p}}Y_j\right), \label{cumulant}
\end{equation}
where the summation is over all partitions
$\bn=(\nu_1,\ldots,\nu_p)$ of size $p=1,\ldots,r$, of the index set
$\{1,2,\ldots,r\}$.

\bigskip
{\sc Assumption 5.} The fourth moment
$E\left[\left(X_{a_j}(t)\right)^4\right]$ is finite for all
$a_j\in\{1,\ldots,r\}$, while the fourth order cumulant function
defined by
$$cum\left[X_{a_1}(t+t_1),X_{a_2}(t+t_2),X_{a_3}(t+t_3),X_{a_4}(t)\right]$$
does not depend on $t$, and this function, denoted by
$C_{a_1a_2a_3a_4}(t_1,t_2,t_3)$, satisfies
$$|C_{a_1a_2a_3a_4}(t_1,t_2,t_3)|\le \prod_{i=1}^{3}g_{a_{i}}(t_i),$$
where $g_{a_{i}}(t_i),\ i=1,2,3,$ are all continuous, even,
nonnegative and integrable functions over the real line, which are
non-increasing over $[0,\infty)$ for all
$a_1,a_2,a_3,a_4\in\{1,2,\ldots,r\}$. %
\bigskip

Note that the cross spectral density is, in general, complex valued.
Thus, the proposed estimator $\widehat{\phi}_{a_1a_2}(\cdot)$ can be
represented as the vector
\begin{equation}\label{est1}
\left( \begin{array}{c}
Re\left(\widehat{\phi}_{a_1a_2}(\lambda)\right)\\[1ex]
Im\left(\widehat{\phi}_{a_1a_2}(\lambda)\right) \end{array} \right),
\end{equation}
where
\begin{eqnarray*}
&&\hskip-20pt Re(\widehat{\phi}_{a_1a_2}(\lambda))\\
&=&\frac{1}{2\pi
n\rho_n}\sum_{t_1=1}^{n}\sum_{t_2=1}^{n}K(b_n(t_2-t_1))X_{a_1}
\left(\frac{t_1}{\rho_n}\right)X_{a_2}\left(\frac{t_2}{\rho_n}\right)\cos\left(\frac{(t_2-t_1)\lambda}{\rho_n}\right)
1_{[-\pi\rho_n,\pi\rho_n]}(\lambda),\\
&&\hskip-20pt Im(\widehat{\phi}_{a_1a_2}(\lambda))\\
&=&\frac{1}{2\pi n\rho_n}\sum_{t_1=1}^{n}\sum_{t_2=1}^{n}K(b_n(t_2-t_1))X_{a_1}%
\left(\frac{t_1}{\rho_n}\right)X_{a_2}\left(\frac{t_2}{\rho_n}\right)\sin\left(\frac{(t_2-t_1)\lambda}{\rho_n}\right)
1_{[-\pi\rho_n,\pi\rho_n]}(\lambda).\\
\end{eqnarray*}

\bigskip{\sc Theorem~2.} {\em Under Assumptions 1--5, the covariance of
$\left(\!\!\!\!\begin{array}{c}
Re\left(\widehat{\phi}_{a_1a_2}(\cdot)\right)\\[1ex]
Im\left(\widehat{\phi}_{a_1a_2}(\cdot)\right)
\end{array}\!\!\!\right)$ %
with %
$\left(\!\!\!\!\begin{array}{c}
Re\left(\widehat{\phi}_{a_3a_4}(\cdot)\right)\\[1ex]
Im\left(\widehat{\phi}_{a_3a_4}(\cdot)\right)
\end{array}\!\!\!\right)$ %
converges as follows:
$$\lim_{n\rightarrow\infty}nb_n Cov\left[\left( \begin{array}{c}
Re\left(\widehat{\phi}_{a_1a_2}(\lambda_1)\right)\\[1ex]
Im\left(\widehat{\phi}_{a_1a_2}(\lambda_1)\right) \end{array}
\right),\left( \begin{array}{c}
Re\left(\widehat{\phi}_{a_3a_4}(\lambda_2)\right)\\[1ex]
Im\left(\widehat{\phi}_{a_3a_4}(\lambda_2)\right) \end{array}
\right)\right] =\left[
\begin{array}{cc}
\sigma_{11}(\lambda_1,\lambda_2)\ \ &
\sigma_{12}(\lambda_1,\lambda_2)\\[1ex]
\sigma_{21}(\lambda_1,\lambda_2)\ \ &
\sigma_{22}(\lambda_1,\lambda_2)
\end{array} \right],
$$%
where
\begin{eqnarray*}
\sigma_{11}(\lambda_1,\lambda_2)&=&B\cdot
Re\left\{\phi_{a_1a_3}^*(\lambda_2)\phi_{a_2a_4}(\lambda_2)+\phi_{a_1a_4}^*(\lambda_2)\phi_{a_2a_3}(\lambda_2)\right\}\\
&&\times[1_{E_2}(\lambda_1,\lambda_2)+1_{E_3}(\lambda_1,\lambda_2)+2\times1_{E_4}(\lambda_1,\lambda_2)],\\
\sigma_{12}(\lambda_1,\lambda_2)&=&B\cdot
Im\left\{\phi_{a_1a_3}(\lambda_2)\phi_{a_2a_4}^*(\lambda_2)+\phi_{a_1a_4}^*(\lambda_2)\phi_{a_2a_3}(\lambda_2)\right\}\\
&&\times[1_{E_2}(\lambda_1,\lambda_2)+1_{E_3}(\lambda_1,\lambda_2)],\\
\sigma_{21}(\lambda_1,\lambda_2)&=&B\cdot
Im\left\{\phi_{a_1a_3}^*(\lambda_2)\phi_{a_2a_4}(\lambda_2)+\phi_{a_1a_4}^{*}(\lambda_2)\phi_{a_2a_3}(\lambda_2)\right\}\\
&&\times[1_{E_2}(\lambda_1,\lambda_2)-1_{E_3}(\lambda_1,\lambda_2)],\\
\sigma_{22}(\lambda_1,\lambda_2)&=&B\cdot
Re\left\{\phi_{a_1a_3}(\lambda_2)\phi_{a_2a_4}^*(\lambda_2)-\phi_{a_1a_4}(\lambda_2)\phi_{a_2a_3}^*(\lambda_2)\right\}\\
&&\times[1_{E_2}(\lambda_1,\lambda_2)-1_{E_3}(\lambda_1,\lambda_2)],\\
B&=&\frac12 \int_{-\infty}^{\infty}K^{2}(x)dx,\\
E_1&=&\{(\lambda_1,\lambda_2):~\lambda_1-\lambda_2\ne0,~\lambda_1+\lambda_2\ne0,~-\infty<\lambda_1,\lambda_2<\infty\},\\
E_2&=&\{(\lambda_1,\lambda_2):~\lambda_1-\lambda_2=0,~-\infty<\lambda_1,\lambda_2<\infty\}~\backslash~\{(0,0)\},\\
E_3&=&\{(\lambda_1,\lambda_2):~\lambda_1+\lambda_2=0,~-\infty<\lambda_1,\lambda_2<\infty\}~\backslash~\{(0,0)\},\\
E_4&=&\{(0,0)\}.
\end{eqnarray*}
The convergence is uniform over any compact subset of $E_1$, $E_2$
or $E_3$. In particular, the variance-covariance matrix of the
random vector %
$\left(\begin{array}{c}
Re\left(\widehat{\phi}_{a_1a_2}(\cdot)\right)\\[1ex]
Im\left(\widehat{\phi}_{a_1a_2}(\cdot)\right) \end{array}
\right)$ %
goes to zero as $n\rightarrow\infty$, for all
$a_1,a_2\in\{1,2,\ldots,r\}$.}
\bigskip

The covariance between two complex-valued random variables is often
defined as the trace of the $2\times2$ cross-covariance matrix of
the random vectors formed by their real and imaginary parts
(Brockwell and Davis, 1991). In the case of the pair
$(\widehat{\phi}_{a_1a_2}(\lambda_1),\widehat{\phi}_{a_3a_4}(\lambda_2))$,
the limiting covariance according to this notion can be easily be
computed from Theorem~2.

Theorem~1 and Theorem~2 together establish the consistency of any
vector of estimators having elements of the form
$\widehat{\phi}_{a_1a_2}(\cdot)$.

\section{Asymptotic Normality}

We will make an additional assumption about the underlying process
in order to prove the asymptotic normality of the estimator.

\bigskip
{\sc Assumption 5A.} The process $\bf{X}$ is strictly stationary;
all moments of the process exist, i.e.,
$E\left[\left(X_{a}(t)\right)^k\right]<\infty$ for each $k>2$ and
for all $a\in\{1,\ldots,r\}$; %
and for each $a_1,a_2,\ldots,a_k\in\{1,2,\ldots,r\}$ and each $k>2$,
the $k$th order joint cumulant denoted by
$$C_{a_1a_2\ldots a_k}(t_1,t_2,\ldots,t_{k-1})
=cum\left(X_{a_1}(t_1+t),X_{a_2}(t_2+t),\ldots,X_{a_{k-1}}(t_{k-1}+t),X_{a_k}(t)\right),$$
satisfies
$$\left|C_{a_1a_2\ldots a_k}(t_1,t_2,\ldots,\ldots,t_{k-1})\right|\le \prod_{i=1}^{k-1}g_{a_{i}}(t_i),$$
where $g_{a_{i}}(t_i),\ i=1,\ldots,k-1$ are continuous, even,
nonnegative and integrable functions over the real line, which are
non-increasing over $(0,\infty)$. %
\bigskip

Note that Assumption~5A is stronger than Assumption~5.

The following theorem describes the asymptotic behaviour of the
joint cumulants of the estimators $\widehat{\phi}_{a_1a_2}(\cdot)$
for $a_1,a_2\in\{1,2,\ldots,r\}$. In the present case, a cumulant
defined as in (\ref{cumulant}) may be complex-valued.

\bigskip{\sc Theorem~3.} {\em Under the Assumptions 1--4 and~5A, for
$L>2$, the $L$th order joint cumulant of the vector
$\left(\widehat{\phi}_{a_1a_2}(\lambda_1),\ldots,\widehat{\phi}_{a_{2L-1}a_{2L}}(\lambda_L)\right)$
for $a_1,a_2,\ldots,a_{2L}\in\{1,2,\ldots,r\}$ is bounded from above
as follows.
\begin{equation}\label{cumrate}
\left|cum\left(\widehat{\phi}_{a_1a_2}(\lambda_1),\ldots,\widehat{\phi}_{a_{2L-1}a_{2L}}(\lambda_L)\right)\right|
\le Q\cdot(nb_n)^{-(L-1)},
\end{equation}
where the constant $Q$ does not depend on
$\lambda_1,\ldots,\lambda_L$.} \bigskip

\bigskip
{\sc Theorem~4}. {\em Under Assumptions 1--4 and~5A, a vector of
real and imaginary parts of estimated spectra or cross-spectra
converges weakly as follows.
\begin{equation}\label{meancorrected}
\sqrt{nb_n}\left[\left(\begin{array}{c}
Re\{\widehat{\phi}_{a_1a_2}(\lambda_1)\}\\
Im\{\widehat{\phi}_{a_1a_2}(\lambda_1)\}\\ \vdots\\
Re\{\widehat{\phi}_{a_{2J-1}a_{2J}}(\lambda_J)\}\\
Im\{\widehat{\phi}_{a_{2J-1}a_{2J}}(\lambda_J)\}\\
\end{array}\right)-E\left(\begin{array}{c}
Re\{\widehat{\phi}_{a_1a_2}(\lambda_1)\}\\
Im\{\widehat{\phi}_{a_1a_2}(\lambda_1)\}\\ \vdots\\
Re\{\widehat{\phi}_{a_{2J-1}a_{2J}}(\lambda_J)\}\\
Im\{\widehat{\phi}_{a_{2J-1}a_{2J}}(\lambda_J)\}\\
\end{array}\right)\right]
\stackrel{D}{\rightarrow}N_{2J}(0,\Sigma),
\end{equation}
where $a_1,a_2,\ldots,a_{2J}\in\{1,2,\ldots,r\}$, and the elements
of $\Sigma$ are defined in accordance with Theorem~2.}
\bigskip

The foregoing theorem only shows that the vector estimator, after
appropriate mean adjustment and scaling, converges weakly to a
multivariate normal distribution. However, weak convergence around
the true vector of spectra and cross-spectra remains to be
established. Note that
\begin{equation}\label{asym1}
\begin{split}
\sqrt{nb_n}\left(\widehat{\phi}_{a_1a_2}(\lambda)-\phi_{a_1a_2}(\lambda)\right)
=&\sqrt{nb_n}\left(\widehat{\phi}_{a_1a_2}(\lambda)-E[\widehat{\phi}_{a_1a_2}(\lambda)]\right)\\
&+\sqrt{nb_n}\left(E[\widehat{\phi}_{a_1a_2}(\lambda)]-\phi_{a_1a_2}(\lambda)\right).
\end{split}
\end{equation}
We make some further assumptions on the smoothness and the rate of
decay of the spectrum and the shape of the kernel function in order
to obtain the rate of convergence of the bias
$E[{\widehat{\phi}_{a_1a_2}}(\lambda)]-{\phi_{a_1a_2}}(\lambda)$.

\bigskip
{\sc Assumption 1A.} The function $g_{qa_1a_2}(\cdot)$, defined over
the real line as $$g_{qa_1a_2}(t)=\sup_{|s|\ge
|t|}|s|^{q}|C_{a_1a_2}(s)|$$ is integrable for all
$a_1,a_2\in\{1,2,\ldots,r\}$, for some positive number $q$ greater
than 1.%

\bigskip
{\sc Assumption 1B.} The power spectral density is such that, for
all $a_1,a_2\in\{1,2,\ldots,r\}$ and for some $p>1$,
$\displaystyle\lim_{\lambda\rightarrow\infty}|\lambda^p\phi_{a_1a_2}(\lambda)|=A_{a_1a_2}$
for some non-negative number $A_{a_1a_2}$.%
\bigskip

For any kernel $K(\cdot)$, let us define
$$k_s=\lim_{x\rightarrow 0}\frac{1-K(x)}{|x|^s}$$
for each positive number $s$ such that the limit exists. The
characteristic exponent of the kernel is defined as the largest
number $s$, such that the limit exists and is non-zero (Parzen,
1957). In other words, the characteristic exponent is the number $s$
such that $1\!-\!K(1/y)$ is $O(y^{-s})$.

\bigskip
{\sc Assumption 2A.} The characteristic exponent of the kernel
$K(\cdot)$ is a number, for which Assumption 1A holds.%
\bigskip

Note that Assumption 1A implies Assumption 1, and also that
$\phi_{a_1a_2}(\cdot)$ is $[q]$ times differentiable, where $[q]$ is
the integer part of $q$. Thus, the number $q$ indicates the degree
of smoothness of the spectral density. If Assumption 1A holds for a
particular value of $q$, then it would also hold for smaller values.

The number $p$ indicates the slowest rate of decay of the various
elements of the power spectral density matrix. The following are two
interesting situations, where Assumption~1B holds.
\begin{enumerate}
\item The real and imaginary parts of the components of the power
spectral density matrix are rational functions of the form
$\frac{P(\lambda)}{Q(\lambda)}$, where $P(\cdot)$ and $Q(\cdot)$ are
polynomials such that the degree of $Q(\cdot)$ is more than degree
of $P(\cdot)$ by at least $p$. Note that continuous time ARMA
processes possess rational power spectral density.
\item The function $C_{a_1a_2}(\cdot)$ has the following smoothness property:
$C_{a_1a_2}(\cdot)$ is $p$ times differentiable and the $p^{\rm th}$
derivative of $C_{a_1a_2}(\cdot)$ is in $L^1$.
\end{enumerate}

\bigskip{\sc Theorem~5.} {\em Under Assumptions 2--4, 1A, 1B and 2A, the
bias of the estimator $\widehat{\phi}_{a_1a_2}(\lambda)$ given by
(\ref{est}), for $a_1,a_2\in\{1,2,\ldots,r\}$, is
\begin{equation*}
\begin{split}
E[\widehat{\phi}_{a_1a_2}(\lambda)-\phi_{a_1a_2}(\lambda)]%
&=\left[-\frac{k_{q}}{2\pi}\int_{-\infty}^{\infty}|t|^{q}C_{a_1a_2}(t)e^{-it\lambda}
dt\right](\rho_n b_n)^{q}+o\left((\rho_n b_n)^{q}\right)\\
&~+\left[-\frac{1}{2\pi}\int_{-\infty}^{\infty}|t|C_{a_1a_2}(t)e^{-it\lambda}
dt\right]\left(\frac{\rho_n}{n}\right)+o\left(\frac{\rho_n}{n}\right)\\
&~+\left[\frac{A_{a_1a_2}}{(2\pi)^{p}}\sum_{|l|>0}\frac{1}{|l|^p}\right]\frac{1}{(\rho_n)^{p}}
+o\left(\frac{1}{(\rho_n)^{p}}\right).
\end{split}
\end{equation*}}
\bigskip

Theorem~5 shows that the second term in (\ref{asym1}) would go to
zero if the sampling rate $\rho_n$ satisfies additional conditions.

\bigskip
{\sc Assumption 4A.} The sampling rate is such that
$\sqrt{nb_n}(\rho_nb_n)^{q}\rightarrow0$ and
$\sqrt{nb_n}/\rho_n^{p}\rightarrow0$ as $n\rightarrow\infty$.%
\bigskip

Note that, whenever Assumption 3 holds, Assumption 4A is stronger
than Assumption~4. With this assumption, the expected values of the
estimators in Theorem~4 can be replaced by the respective true
values.

\bigskip{\sc Theorem~6.} {\em Under Assumptions 1--3, 1A, 1B, 2A, 4A and
5A, we have the following weak convergence.
$$\sqrt{nb_n}\left[\left(\begin{array}{c}
Re\{\widehat{\phi}_{a_1a_2}(\lambda_1)\}\\
Im\{\widehat{\phi}_{a_1a_2}(\lambda_1)\}\\ \vdots\\
Re\{\widehat{\phi}_{a_{2J-1}a_{2J}}(\lambda_J)\}\\
Im\{\widehat{\phi}_{a_{2J-1}a_{2J}}(\lambda_J)\}\\
\end{array}\right)-\left(\begin{array}{c}
Re\{{\phi}_{a_1a_2}(\lambda_1)\}\\
Im\{{\phi}_{a_1a_2}(\lambda_1)\}\\ \vdots\\
Re\{{\phi}_{a_{2J-1}a_{2J}}(\lambda_J)\}\\
Im\{{\phi}_{a_{2J-1}a_{2J}}(\lambda_J)\}\\
\end{array}\right)\right]
\stackrel{D}{\rightarrow}N_{2J}(0,\Sigma),$$ %
where $a_1,a_2,\ldots,a_{2J}\in\{1,2,\ldots,r\}$, and the elements
of $\Sigma$ are defined in accordance with Theorem~2.}

\section{Optimal rate of convergence}
We are now in a position to optimize the rates of $b_n$ and $\rho_n$
so that $\frac{1}{\sqrt{nb_n}}$ tends to 0 as fast as possible under
the conditions of Theorem~6.

\bigskip{\sc Theorem~7.} {\em Under Assumptions 3 and 4A, the reciprocal
of the scale factor ($\frac{1}{\sqrt{nb_n}}$) used in Theorem~6 has
the fastest convergence to 0 when
\begin{eqnarray*}
b_n&=& o\left(n^{-\frac{p+q}{p+q+2pq}}\right),\\
\rho_n&=& O\left(n^{\frac q{p+q+2pq}}\right),
\end{eqnarray*}
and under these conditions, $\frac1{\sqrt{nb_n}}=o\left(n^{-\frac
{pq}{p+q+2pq}}\right)$.}
\bigskip

It has been shown in Srivastava and Sengupta (2010) that under the
assumptions of Theorems~2 and~5, the optimal rate of convergence for
mean square consistency of the estimator (\ref{est}) is given as
$$E\left[\{\widehat\phi_{a_1a_2}(\cdot)-\phi_{a_1a_2}(\cdot)\}^2\right]=O\left(n^{-\frac
{2pq}{p+q+2pq}}\right),$$ %
which corresponds to the choices
\begin{eqnarray*}
b_n&=&O\left(n^{-\frac{p+q}{p+q+2pq}}\right),\\
\rho_n&=& O\left(n^{\frac q{p+q+2pq}}\right).
\end{eqnarray*}
Theorem~7 shows that the optimal rate of weak convergence of the
estimator $\widehat{\phi}_{a_1a_2}(\cdot)$ is slower than the square
root of the optimal rate corresponding to mean square consistency.

It is important to note that for every fixed value of~$q$, the
number $p$, which indicates rate of decay of the spectrum, can be
increased indefinitely by continuous time low pass filtering with a
cut off frequency larger than the maximum frequency of interest.
There are well-known filters such as the Butterworth filter, which
have polynomial rate of decay of the transfer function with
specified degree of the polynomial, that can be used for this
purpose. For fixed $q$, the best rate of weak convergence given in
Theorem~7, obtained by allowing $p$ to go to infinity, happens to be
$o\left(n^{-\frac {q}{1+2q}}\right)$.

The rate of weak convergence crucially depends on the number $q$,
the assumed degree of smoothness of the spectrum. The stronger the
assumption, the faster is the rate of convergence. The rate
corresponding to $q=1$ (weakest possible assumption) is
$o\left(n^{-\frac13}\right)$, assuming that $p$ can be allowed to be
very large. For very large $q$ (very strong assumption) and large
$p$, the rate approaches $o\left(n^{-\frac12}\right)$.

\section{Simulation}
With a view to investigating the applicability of the asymptotic
results reported in Section~3 to finite sample size, we consider the
bivariate continuous time linear process
$$
\left(\begin{array}{c}
X_1(t)\\[2ex]
X_2(t)\\
\end{array}\right)=\left(\begin{array}{c}
\int_{-\infty}^{t}h_{1}(t-u)Z_1(u)du+\int_{-\infty}^{t}h_{2}(t-u)Z_2(u)du\\[2ex]
\int_{-\infty}^{t}h_{3}(t-u)Z_1(u)du+\int_{-\infty}^{t}h_{4}(t-u)Z_3(u)du\\
\end{array}\right)
$$
where $Z_j(u)$, $j=1,2,3$ are independent continuous time white
noise and $h_{j}(u)=\beta_je^{-\alpha_ju}$ for $j\in\{1,2,3,4\}$.
The elements of the spectral density matrix
$$\left(\begin{array}{cc}\phi_{11}(\lambda)&\phi_{12}(\lambda)\\
\phi_{12}^*(\lambda)&\phi_{22}(\lambda)\\
\end{array}\right)$$
are defined as follows (Hoel et al., 1972).
\begin{eqnarray*}
\phi_{11}(\lambda)&=&\frac{1}{2\pi}\cdot\frac{\beta_1^{2}}{\alpha_1^2+\lambda^2}+\frac{1}{2\pi}\cdot\frac{\beta_2^{2}}{\alpha_2^2+\lambda^2},\\
\phi_{22}(\lambda)&=&\frac{1}{2\pi}\cdot\frac{\beta_3^{2}}{\alpha_3^2+\lambda^2}+\frac{1}{2\pi}\cdot\frac{\beta_4^{2}}{\alpha_4^2+\lambda^2},\\
Re(\phi_{12}(\lambda))&=&\frac{1}{2\pi}\cdot\frac{\beta_1\beta_3(\alpha_1\alpha_3+\lambda^2)}{(\alpha_1^2+\lambda^2)(\alpha_3^2+\lambda^2)}\\
\mbox{and
}Im(\phi_{12}(\lambda))&=&\frac{1}{2\pi}\cdot\frac{\beta_1\beta_3(\alpha_3-\alpha_1)\lambda}{(\alpha_1^2+\lambda^2)(\alpha_3^2+\lambda^2)}.
\end{eqnarray*}
We simulate this bivariate process with the choices $\beta_1=1$,
$\beta_2=1$, $\beta_3=2$, $\beta_4=\frac 2 5$,
$\alpha_1=\beta_1\cdot\sqrt{\frac 3 2}$,
$\alpha_2=\beta_2\cdot\sqrt{3}$, $\alpha_3=\beta_3\cdot\sqrt{3}$ and
$\alpha_4=\beta_4\cdot\sqrt{3}$. Note that for this process,
Assumption 1A holds with $q\ge1$ and Assumption 1B holds with
$p\le2$. For the purpose of estimation, we make these assumptions
with $p=2$ and $q=2$. In accordance with this choice of $q$, we use
the second order kernel function
$$K(x)=\frac 1 2 \{1+\cos(\pi x)\}1_{[-1,1]}(x).$$
We also use the rates $b_n=\frac 1 4 n^{-\frac 1 4}$ and
$\rho_n=4\cdot n^{\frac{1}{6}}$.

We estimate the bivariate spectrum matrix for frequencies in the
range $[0,3\pi]$ at intervals of $.01\pi$ (i.e., $301$ uniformly
spaced grid points). We subsequently compute the normalized
statistics
\begin{eqnarray*}
T_1(\lambda)&=&\sqrt{nb_n}\left(\frac{\widehat\phi_{11}(\lambda)-\phi_{11}(\lambda)}
{\sqrt{2\{1+1_{\{0\}}(\lambda)\}B\widehat\phi_{11}^2(\lambda)}}\right),\\
T_2(\lambda)&=&\sqrt{nb_n}\left(\frac{\widehat\phi_{22}(\lambda)-\phi_{22}(\lambda)}
{\sqrt{2\{1+1_{\{0\}}(\lambda)\}B\widehat\phi_{22}^2(\lambda)}}\right),\\
T_3(\lambda)&=&\sqrt{nb_n}\left(\frac{Re(\widehat\phi_{12}(\lambda))-Re(\phi_{12}(\lambda))}
{\sqrt{\{1+1_{\{0\}}(\lambda)\}B[\widehat\phi_{11}(\lambda)\widehat\phi_{22}(\lambda)+\{Re(\widehat\phi_{12}(\lambda))\}^2
-\{Im(\widehat\phi_{12}(\lambda))\}^2]}}\right),\\
T_4(\lambda)&=&\sqrt{nb_n}\left(\frac{Im(\widehat\phi_{12}(\lambda))-Im(\phi_{12}(\lambda))}
{\sqrt{B[\widehat\phi_{11}(\lambda)\widehat\phi_{22}(\lambda)-\{Re(\widehat\phi_{12}(\lambda))\}^2
+\{Im(\widehat\phi_{12}(\lambda))\}^2]}}\right)\left[1-1_{\{0\}}(\lambda)\right],\\
\end{eqnarray*}
in accordance with Theorem~2. According to Theorem~6, the asymptotic
distribution of each of these four statistics is standard normal.
This procedure is repeated for 500 simulation runs. By regrading the
values of the above statistics for the different simulation runs as
four data sets of size 500 each, we calculate the Kolmogorov-Smirnov
test statistic (Shorak and Wellner, 1986) for these data sets, and
the corresponding p-value. This procedure is repeated for the 301
frequency values mentioned above. The percentage of p-values (across
301 frequency values) exceeding the number 0.05 are reported in
Table~1, for sample sizes $n=100$, 1000, 10000 and 100000. The table
shows that for each statistic, the percentage approaches the ideal
value of 95 very slowly as $n$ increases.

We now turn to computation of confidence limits of the power
spectral density. For each frequency value, we compute the 95\%
asymptotic confidence intervals of $\phi_{11}$, $\phi_{22}$,
$Re(\phi_{12})$ and $Im(\phi_{12})$ from the statistics
$T_1(\lambda)$, $T_2(\lambda)$, $T_3(\lambda)$ and $T_4(\lambda)$,
assuming that the latter have the standard normal distribution.
Subsequently, we compute the fraction of times (out of 500
simulation runs) the confidence limits capture the true value of the
function. These percentages are plotted against the frequency, for
sample sizes $n=100$, 1000, 10000 and 100000, in Figure~1. It is
seen that the observed fraction approaches the ideal coverage
probability (0.95) for larger sample sizes.
 Since there is a discontinuity of the asymptotic variance
function at the point $\lambda=0$, while the estimated spectrum is
constrained to be continuous, some anomalous behaviour\ \,in the
neighbourhood

\begin{center}
\begin{tabular}{|c|rrrr|}
\hline sample& \multicolumn{4}{|c|}{observed percentage}\\
\cline{2-5}
size ($n$) &$\phi_{11}$&$\phi_{22}$&$Re(\phi_{12})$&$Im(\phi_{12})$\\
\hline
100 &   0.0 \%  &   0.0 \%  &   0.0 \%  &   0.0 \%  \\
1000    &   4.3 \%  &   6.0 \%  &   9.6 \%  &   16.3    \%  \\
10000   &   73.4    \%  &   72.4    \%  &   76.1    \%  &   72.8    \%  \\
100000  &   91.4    \%  &   88.4    \%  &   93.4    \%  &   90.0    \%  \\
\hline
\end{tabular}
\end{center}
\noindent{\small Table 1. Observed percentage of frequencies (in the
range 0 to $3\pi$) for which p-values of the Kolmogorov-Smirnov
statistics for testing normality of $\phi_{11}$, $\phi_{22}$,
$Re(\phi_{12})$ and $Im(\phi_{12})$ are greater than 0.05 (ideal
percentage is 95\%).}

\noindent\includegraphics[width=6in]{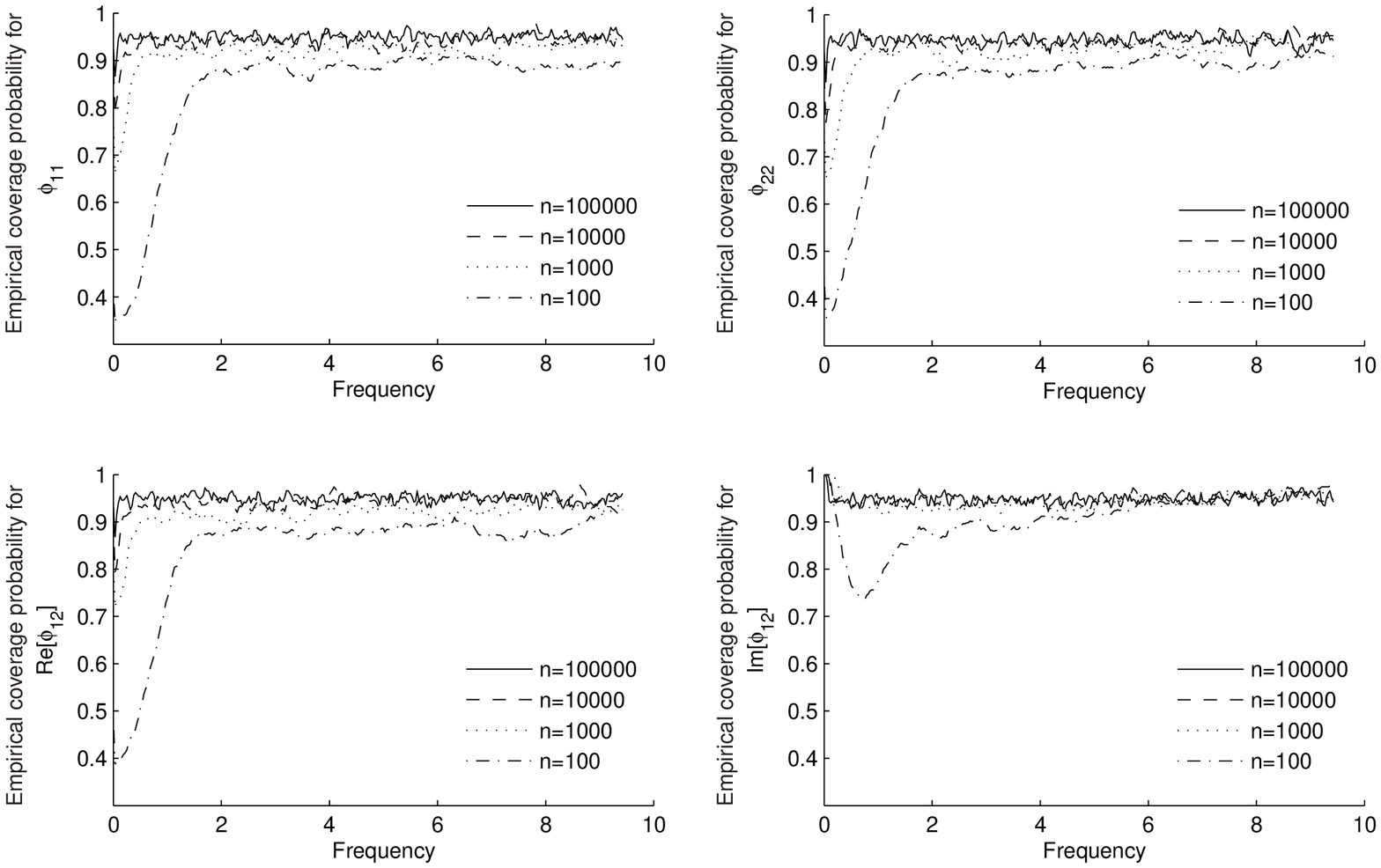}

\noindent {\small Figure 1. Empirical coverage probability (based on
500 simulation runs) of pointwise confidence intervals of
$\phi_{11}$, $\phi_{22}$, $Re(\phi_{12})$ and $Im(\phi_{12})$ for
sample sizes 100, 1000, 10000 and 100000.}

\bigskip
\noindent of the the point $\lambda=0$ is expected. This results in
substantially lower values of the empirical coverage probability in
this region. However, this region of anomaly is observed to shrink
as the sample size increases. It would be interesting to note that
the empirical coverage probability is reasonably close to the ideal
coverage probability for most frequency values when the sample size
as small as 1000, even though Table~1 indicates that the asymptotic
distribution is not applicable at this sample size.

\section*{Appendix} \setcounter{section}{1}\setcounter{equation}{0}
\renewcommand{\theequation}{\Alph{section}.\arabic{equation}}

We denote by $K_1(\cdot)$ a function that bounds the covariance
averaging kernel $K(\cdot)$ as in Assumption~2. Further, we denote
$K_1(0)$ by $M$.

\medskip{\sc Proof of Theorem~1.} We shall show that the bias of the
estimator $\widehat{\phi}_{a_1a_2}(\lambda)$ given by (\ref{est})
converges to $0$ uniformly over $[\lambda_l,\lambda_u]$ for any
$\lambda_l$, $\lambda_u$ such that $\lambda_l<\lambda_u$. Note that
\begin{equation*}
E[\widehat{\phi}_{a_1a_2}(\lambda)] =\frac{1}{2\pi
\rho_n}\sum_{u=-(n-1)}^{n-1}\left(1-\frac{|u|}{\rho_n}\right)K(b_nu)C_{a_1a_2}
\left(\frac{u}{\rho_n}\right)e^{-\frac{iu\lambda}{\rho_n}}1_{[-\pi\rho_n,\pi\rho_n]}(\lambda).
\end{equation*}
Consider the simple function $S_n(\cdot)$, defined over
$[\lambda_l,\lambda_u]\times\mathbb{R}$, by
\begin{align*}
S_n(\lambda,x)=\frac{1}{2\pi
}\sum_{u=-(n-1)}^{n-1}\left(1-\frac{|u|}{\rho_n}\right)K(b_nu)C_{a_1a_2}
\left(\frac{u}{\rho_n}\right)e^{-\frac{iu\lambda}{\rho_n}}1_{[-\pi\rho_n,\pi\rho_n]}(\lambda)1_{\left(\frac{u-1}{\rho_n},\frac{u}{\rho_n}\right]}(x).
\end{align*}
Observe that $\int_{-\infty}^{\infty}
S_n(\lambda,x)dx=E[\widehat{\phi}_{a_1a_2}(\lambda)]$. Define the
function $S(\cdot)$, over $[\lambda_l,\lambda_u]\times\mathbb{R}$,
by $ S(\lambda,x)=\frac{1}{2\pi}C_{a_1a_2}(x)e^{-ix\lambda}.$

For any $x \in \mathbb{R}$, let $u_{n}(x)$ be the smallest integer
greater than or equal to $\rho_nx$. Note that the interval
$\left(\frac{u_{n-1}(x)}{\rho_n},\frac{u_{n}(x)}{\rho_n}\right]$
contains the point $x$  and
$\lim_{n\rightarrow\infty}\frac{u_{n}(x)}{\rho_n}= x$. For
sufficiently large $n$, we have from Assumptions 3 and 4,
\begin{align*}
S_n(\lambda,x)=&\frac{1}{2\pi}\left(1-\frac{|u_{n}(x)|}{\rho_n}\frac{\rho_n}{n}\right)K\left(b_n
\rho_n\frac{u_{n}(x)}{\rho_n}\right)
C_{a_1a_2}\left(\frac{u_{n}(x)}{\rho_n}\right)
e^{-\frac{iu_{n}(x)\lambda}{\rho_n}}1_{[-\pi\rho_n,\pi\rho_n]}(\lambda).
\end{align*}
Proving the uniform convergence of
$Bias[\widehat{\phi}_{a_1a_2}(\lambda)]$ over the finite interval
$[\lambda_l,\,\lambda_u]$ amounts to proving
\begin{equation}\label{alt}
\lim_{n\rightarrow\infty}\int_{-\infty}^{\infty}
S_n(\lambda,x)dx=\int_{-\infty}^{\infty} S(\lambda,x)dx,
\end{equation}
uniformly over  $[\lambda_l,\,\lambda_u]$.

Observe that $\int_{-\infty}^{\infty}
S(\lambda,t)dt=\phi_{a_1a_2}(\lambda)$, which is continuous. By
virtue of the continuity of the limiting function, (\ref{alt}) is
equivalent to proving that $\int_{-\infty}^{\infty}
S_n(\lambda,x)dx$ converges continuously over this interval
(Resnick, 1987), i.e., for any sequence $\lambda_n \rightarrow
\lambda$,
\begin{equation}
\lim_{n\rightarrow\infty}\int_{-\infty}^{\infty} S_n(\lambda_n,x)dx
= \int_{-\infty}^{\infty} S(\lambda,x)dx, \label{star}
\end{equation}
where $\lambda_n,\lambda \in  [\lambda_l,\,\lambda_u]$.

By continuity of the function $S_n(\lambda,x)$ with respect to $x$
and $\lambda$, we have from Assumptions 3 and 4, for any fixed $x$,
$$\lim_{n\rightarrow \infty}|S_n(\lambda_n,x)-S(\lambda,x)|=0.$$
Note that from Assumptions 1 and 2, we have the dominance
\begin{align*}
|S_n(\lambda_n,x)|\le&
M\sum_{|u|<n}\left|C_{a_1a_2}\left(\frac{u}{\rho_n}\right)\right|1_{\left(\frac{u-1}{\rho_n},\frac{u}{\rho_n}\right]}(x)\le
M g_{a_1a_2}(x),
\end{align*}
where $g_{a_1a_2}(\cdot)$ is the function described in Assumption~1.
Thus, by applying the dominated convergence theorem (DCT), we have
(\ref{star}).

Hence, $E[\widehat{\phi}_{a_1a_2}(\lambda)]\rightarrow
\phi(\lambda)$ uniformly on $[\lambda_l,\,\lambda_u].\hfill\Box$

\medskip{\sc Proof of Theorem~2.} We begin by calculating the covariance
between the estimators $Re(\widehat{\phi}_{a_1a_2}(\cdot))$ and
$Re(\widehat{\phi}_{a_3a_4}(\cdot))$.
\begin{equation*}
\begin{split}
&\hskip-8pt Cov\left[Re(\widehat{\phi}_{a_1a_2}(\lambda_1)),Re(\widehat{\phi}_{a_3a_4}(\lambda_2))\right]\\
=&\frac{1}{(2\pi)^{2}(n\rho_n)^{2}}
\sum_{t_1=1}^{n}\sum_{t_2=1}^{n}\sum_{t_3=1}^{n}\sum_{t_4=1}^{n}K(b_n(t_2-t_1))K(b_n(t_4-t_3))\\&\times
Cov\left[X_{a_1}\left(\frac{t_1}{\rho_n}\right)X_{a_2}\left(\frac{t_2}{\rho_n}\right),
X_{a_3}\left(\frac{t_3}{\rho_n}\right)X_{a_4}\left(\frac{t_4}{\rho_n}\right)\right]
\cos\left(\frac{(t_2-t_1)\lambda_1}{\rho_n}\right)\cos\left(\frac{(t_4-t_3)\lambda_2}{\rho_n}\right)\\
=&\frac{1}{(2\pi)^{2}(n\rho_n)^{2}}
\sum_{t_1=1}^{n}\sum_{t_2=1}^{n}\sum_{t_3=1}^{n}\sum_{t_4=1}^{n}K(b_n(t_1-t_2))K(b_n(t_3-t_4))\\&~\times
\left[C_{a_1a_3}\left(\frac{t_1-t_3}{\rho_n}\right)
C_{a_2a_4}\left(\frac{t_2-t_4}{\rho_n}\right)
+C_{a_1a_4}\left(\frac{t_1-t_4}{\rho_n}\right)
C_{a_2a_3}\left(\frac{t_2-t_3}{\rho_n}\right)\right.\\&~~~\left.
+C_{a_1a_2a_3a_4}\left(\frac{t_1-t_4}{\rho_n},\frac{t_2-t_4}{\rho_n},\frac{t_3-t_4}{\rho_n}\right)\right]
\cos\left(\frac{(t_1-t_2)\lambda_1}{\rho_n}\right)\cos\left(\frac{(t_3-t_4)\lambda_2}{\rho_n}\right)
\\
=&T_1(\lambda_1,\lambda_2)+T_2(\lambda_1,\lambda_2)+T_3(\lambda_1,\lambda_2),
\end{split}
\end{equation*}
where the three terms correspond to the three summands appearing
inside square brackets in the previous step.

Now consider the function $T_1(\lambda_1,\lambda_2)$. By using the
transformations $u_1=t_1-t_2$, $u_2=t_1-t_3$ and $u_3=t_2-t_4$, we
have
\begin{equation*}
\begin{split}
T_1(\lambda_1,\lambda_2)=&\frac{1}{(2\pi)^{2}(n\rho_n)^{2}}\sum_{t_1=1}^{n}\sum_{u_1=t_1-1}^{t_1-n}\sum_{u_2=t_1-1}^{n-t_1}\sum_{u_3=t_1-1-u_1}^{t_1-n-u_1}K(b_n
u_1)K(b_n(u_1-u_2+u_3))\\& ~\times
C_{a_1a_3}\left(\frac{u_2}{\rho_n}\right)C_{a_2a_4}\left(\frac{u_3}{\rho_n}\right)
\cos\left(\frac{u_1\lambda_1}{\rho_n}\right)\cos\left(\frac{(u_1-u_2+u_3)\lambda_2}{\rho_n}\right).
\end{split}
\end{equation*}
The range of the four summations on the right hand side is described
by the set of inequalities $1\le t_1\le n$ and $t_1-n\le u_1, u_2,
u_1+u_3\le t_1-1$, which is equivalent to the inequalities
$-(n-1)\le u_1, u_2, u_1+u_3\le (n-1)$ and $\max\{u_1, u_2,
u_1+u_3\}+1\le t_1 \le \min\{u_1, u_2, u_1+u_3\}$. Therefore, the
expression for $T_1(\lambda_1,\lambda_2)$ simplifies to
\begin{equation*}
\begin{split}
&\frac{1}{(2\pi)^{2}n\rho_n^{2}}\sum_{u_1=-(n-1)}^{(n-1)}\sum_{u_2=-(n-1)}^{(n-1)}\sum_{u_3=-(n-1)-u_1}^{(n-1)-u_1}U_n(u_1,u_2,u_3)K(b_n
u_1)K(b_n(u_1-u_2+u_3))\\& ~\times
C_{a_1a_3}\left(\frac{u_2}{\rho_n}\right)C_{a_2a_4}\left(\frac{u_3}{\rho_n}\right)
\cos\left(\frac{u_1\lambda_1}{\rho_n}\right)\cos\left(\frac{(u_1-u_2+u_3)\lambda_2}{\rho_n}\right),
\end{split}
\end{equation*}
where
$$U_n(u_1,u_2,u_3)=\left(1+\frac{\min(u_1,u_2,u_1+u_3)}{n}-\frac{\max(u_1,u_2,u_1+u_3)}{n}\right).$$
By writing the cosine functions in terms of complex exponentials, we
have
\begin{equation}
\begin{split}
&\hskip-8pt T_{1}(\lambda_1,\lambda_2)\\
=&\frac{1}{(4\pi)^{2}n\rho_n^{2}}\sum_{u_1=-(n-1)}^{(n-1)}\sum_{u_2=-(n-1)}^{(n-1)}\sum_{u_3=-(n-1)-u_1}^{(n-1)-u_1}U_n(u_1,u_2,u_3)K(b_n
u_1)K(b_n(u_1-u_2+u_3))\\&~\times
C_{a_1a_3}\left(\frac{u_2}{\rho_n}\right)
C_{a_2a_4}\left(\frac{u_3}{\rho_n}\right)
\left\{e^{-i\frac{(\lambda_1-\lambda_2)u_1}{\rho_n}}e^{-i\frac{\lambda_2u_2}{\rho_n}}e^{i\frac{\lambda_2u_3}{\rho_n}}
+e^{i\frac{(\lambda_1-\lambda_2)u_1}{\rho_n}}e^{i\frac{\lambda_2u_2}{\rho_n}}e^{-i\frac{\lambda_2u_3}{\rho_n}}\right.\\&
~~~\left.+e^{i\frac{(\lambda_1+\lambda_2)u_1}{\rho_n}}e^{-i\frac{\lambda_2u_2}{\rho_n}}e^{i\frac{\lambda_2u_3}{\rho_n}}+e^{-i\frac{(\lambda_1+\lambda_2)u_1}{\rho_n}}e^{i\frac{\lambda_2u_2}{\rho_n}}e^{-i\frac{\lambda_2u_3}{\rho_n}}\right\}\\&
=T_{11}(\lambda_1,\lambda_2)+T_{12}(\lambda_1,\lambda_2)+T_{13}(\lambda_1,\lambda_2)+T_{14}(\lambda_1,\lambda_2),
\end{split}
\label{eq54}
\end{equation}
where the four terms correspond to the four summands appearing
within braces in the last factor on the right hand side of
(\ref{eq54}).

By using the results of Lemmas 1 and~2 given below, we have the
convergence
$$\lim_{n\rightarrow\infty}nb_nT_{11}(\lambda_1,\lambda_2)=
\frac14\left(\int_{-\infty}^{\infty}K^{2}(x)dx\right)\phi_{a_1a_3}(\lambda_2)\phi_{a_2a_4}^*(\lambda_2)
1_{E_2\cup E_4}(\lambda_1,\lambda_2)$$%
and similar arguments show that
\begin{eqnarray*}
\lim_{n\rightarrow\infty} nb_nT_{12}(\lambda_1,\lambda_2)&=&
\frac14\left(\int_{-\infty}^{\infty}K^{2}(x)dx\right)\phi_{a_1a_3}^*(\lambda_2)\phi_{a_2a_4}(\lambda_2)
1_{E_2\cup E_4}(\lambda_1,\lambda_2),\\
\lim_{n\rightarrow\infty} nb_nT_{13}(\lambda_1,\lambda_2)&=&
\frac14\left(\int_{-\infty}^{\infty}K^{2}(x)dx\right)\phi_{a_1a_3}(\lambda_2)\phi_{a_2a_4}^*(\lambda_2)
1_{E_3\cup E_4}(\lambda_1,\lambda_2)\mbox{ and}\\
\lim_{n\rightarrow\infty}nb_nT_{14}(\lambda_1,\lambda_2)&=&
\frac14\left(\int_{-\infty}^{\infty}K^{2}(x)dx\right)\phi_{a_1a_3}^*(\lambda_2)\phi_{a_2a_4}(\lambda_2)
1_{E_3\cup E_4}(\lambda_1,\lambda_2).
\end{eqnarray*}

For the function $T_2(\lambda_1,\lambda_2)$, one can similarly use
the transformations $u_1=t_1-t_2$, $u_2=t_1-t_4$ and $u_3=t_2-t_3$,
interchange the order of summation and expand the cosine functions
in terms of complex exponentials to obtain
\begin{equation*}
\begin{split}
&\hskip-8pt T_{2}(\lambda_1,\lambda_2)\\
=&\frac{1}{(2\pi)^{2}n\rho_n^{2}}\sum_{u_1=-(n-1)}^{(n-1)}\sum_{u_2=-(n-1)}^{(n-1)}\sum_{u_3=-(n-1)-u_1}^{(n-1)-u_1}U_n(u_1,u_2,u_3)K(b_n
u_1)K(b_n(-u_1+u_2-u_3))\\
&~\times
C_{a_1a_4}\left(\frac{u_2}{\rho_n}\right)C_{a_2a_3}\left(\frac{u_3}{\rho_n}\right)
\cos\left(\frac{u_1\lambda_1}{\rho_n}\right)\cos\left(\frac{(-u_1+u_2-u_3)\lambda_2}{\rho_n}\right)\\
=&\frac{1}{(4\pi)^{2}n\rho_n^{2}}\sum_{u_1=-(n-1)}^{(n-1)}\sum_{u_2=-(n-1)}^{(n-1)}\sum_{u_3=-(n-1)-u_1}^{(n-1)-u_1}
U_n(u_1,u_2,u_3)K(b_n u_1)K(b_n(-u_1+u_2-u_3))\\& \times
C_{a_1a_4}\left(\frac{u_2}{\rho_n}\right)C_{a_2a_3}\left(\frac{u_3}{\rho_n}\right)
\left\{e^{-i\frac{(\lambda_1-\lambda_2)u_1}{\rho_n}}e^{-i\frac{\lambda_2u_2}{\rho_n}}e^{i\frac{\lambda_2u_3}{\rho_n}}
+e^{i\frac{(\lambda_1-\lambda_2)u_1}{\rho_n}}e^{i\frac{\lambda_2u_2}{\rho_n}}e^{-i\frac{\lambda_2u_3}{\rho_n}}\right.\\
&~~~\left.+e^{i\frac{(\lambda_1+\lambda_2)u_1}{\rho_n}}e^{-i\frac{\lambda_2u_2}{\rho_n}}e^{i\frac{\lambda_2u_3}{\rho_n}}+e^{-i\frac{(\lambda_1+\lambda_2)u_1}{\rho_n}}e^{i\frac{\lambda_2u_2}{\rho_n}}e^{-i\frac{\lambda_2u_3}{\rho_n}}\right\}\\&
=T_{21}(\lambda_1,\lambda_2)+T_{22}(\lambda_1,\lambda_2)+T_{23}(\lambda_1,\lambda_2)+T_{24}(\lambda_1,\lambda_2).
\end{split}
\end{equation*}
By using similar arguments as in the case of
$nb_nT_{11}(\lambda_1,\lambda_2)$, it can be shown that
\begin{eqnarray*}
nb_nT_{21}(\lambda_1,\lambda_2)&=&
\frac14\left(\int_{-\infty}^{\infty}K^{2}(x)dx\right)\phi_{a_1a_4}(\lambda_2)\phi_{a_2a_3}^*(\lambda_2)%
1_{E_2\cup E_4}(\lambda_1,\lambda_2),\\
nb_nT_{22}(\lambda_1,\lambda_2)&=&
\frac14\left(\int_{-\infty}^{\infty}K^{2}(x)dx\right)\phi_{a_1a_4}^*(\lambda_2)\phi_{a_2a_3}(\lambda_2)%
1_{E_2\cup E_4}(\lambda_1,\lambda_2),\\
nb_nT_{23}(\lambda_1,\lambda_2)&=&
\frac14\left(\int_{-\infty}^{\infty}K^{2}(x)dx\right)\phi_{a_1a_4}(\lambda_2)\phi_{a_2a_3}^*(\lambda_2)%
1_{E_3\cup E_4}(\lambda_1,\lambda_2)\mbox{ and}\\
nb_nT_{24}(\lambda_1,\lambda_2)&=&
\frac14\left(\int_{-\infty}^{\infty}K^{2}(x)dx\right)\phi_{a_1a_4}^*(\lambda_2)\phi_{a_2a_3}(\lambda_2)%
1_{E_3\cup E_4}(\lambda_1,\lambda_2).
\end{eqnarray*}

Finally, for the term $T_3(\lambda_1,\lambda_2)$, we use the
transformations $u_1=t_1-t_4$, $u_2=t_2-t_4$ and $u_3=t_3-t_4$ and
interchange the order of summations to have
\begin{equation*}
\begin{split}
T_3(\lambda_1,\lambda_2)&=\frac{1}{(2\pi)^{2}(n\rho_n)^{2}}\sum_{u_1=-(n-1)}^{(n-1)}\sum_{u_2=-(n-1)}^{n-1}\sum_{u_3=-(n-1)}^{(n-1)}
(n-\min(u_1,u_2,u_3)+\max(u_1,u_2,u_3))\\& K(b_n(u_1-u_2))K(b_n u_3)
C_{a_1a_2a_3a_4}\left(\frac{u_1}{\rho_n},\frac{u_2}{\rho_n},\frac{u_3}{\rho_n}\right)
\cos\left(\frac{(u_1-u_2)\lambda_1}{\rho_n}\right)\cos\left(\frac{u_3
\lambda_2}{\rho_n}\right).
\end{split}
\end{equation*}
From Assumptions~2 and~5, we have
\begin{equation}\label{I3bound}
\begin{split}
nb_n|T_3(\lambda_1,\lambda_2)|\le \rho_n b_n
M^{2}\sum_{u_1=-(n-1)}^{n-1}\sum_{u_2=-(n-1)}^{n-1}\sum_{u_2=-(n-1)}^{n-1}
g_{a_1}\left(\frac{u_1}{\rho_n}\right)g_{a_2}\left(\frac{u_2}{\rho_n}\right)g_{a_3}\left(\frac{u_3}{\rho_n}\right)\frac{1}{\rho_n^3}.
\end{split}
\end{equation}
Now consider the function $S_n(\cdot)$ defined over $\mathbb{R}$ as
\begin{equation*}
\begin{split}
S_n(x)=&\sum_{u_1=-(n-1)}^{n-1}g_{a_1}\left(\frac{u_1}{\rho_n}\right)1_{(\frac{u_1-1}{\rho_n},\frac{u_1}{\rho_n}]}(x).
\end{split}
\end{equation*}
Observe that $\lim_{n\rightarrow\infty}S_n(x)=g_{a_1}(x)$ and
$|S_n(x)|$ is dominated by $g_{a_1}(\cdot)$. By applying DCT, we
have
$$\lim_{n\rightarrow\infty}\int_{-\infty}^{\infty}S_n(x)dx
=\lim_{n\rightarrow\infty}\sum_{u_1=-(n-1)}^{n-1}g_{a_1}\left(\frac{u_1}{\rho_n}\right)
\frac{1}{\rho_n}=\int_{-\infty}^{\infty}g_{a_1}(x)dx.$$ %
Thus, the upper bound of $nb_nT_3(\lambda_1,\lambda_2)$ given by
(\ref{I3bound}) is $O(\rho_n b_n)$. Assumption~\ref{As4} ensures
that $nb_nT_3(\lambda_1,\lambda_2)$ converges to zero uniformly.

By combining all these terms, we have the convergence of
$nb_nCov\left[Re(\widehat{\phi}_{a_1a_2}(\lambda_1)),Re(\widehat{\phi}_{a_3a_4}(\lambda_2))\right]$
as given in the theorem. Convergence of the other three covariances
follow from a similar argument.\\\mbox{} \hfill$\Box$

\bigskip{\sc Lemma~1.} {\em For $\lambda_1-\lambda_2=0$, the function
$T_{11}(\lambda_1,\lambda_2)$ converges as follows.
$$\lim_{n\rightarrow\infty} nb_nT_{11}(\lambda_1,\lambda_2)=\frac14\left(\int_{-\infty}^{\infty}K^{2}(x)dx\right)
\phi_{a_1a_3}(\lambda_2)\phi_{a_2a_4}^*(\lambda_2).$$ The
convergence is uniform on any compact subset of the set
$$E=\{(\lambda_1,\lambda_2):\lambda_1-\lambda_2=0,~-\infty<\lambda_1,\lambda_2<\infty\}.$$}

{\sc Proof of Lemma~1.} Consider a compact subset $E'$ of the set $E$.
Consider the simple function $S_n(\cdot)$, defined over
$E'\times\mathbb{R}^3$ by
\begin{equation*}
\begin{split}
&\hskip-8pt S_n(\lambda_1,\lambda_2,x_1,x_2,x_3)\\
=&\sum_{u_1=-(n-1)}^{(n-1)}\sum_{u_2=-(n-1)}^{(n-1)}\sum_{u_3=-(n-1)-u_1}^{(n-1)-u_1}U_n(u_1,u_2,u_3)
K(b_n u_1)K(b_n(u_1-u_2+u_3))
\\&~\times  C_{a_1a_3}\left(\frac{u_2}{\rho_n}\right)e^{-i\frac{u_2\lambda_2}{\rho_n}} C_{a_2a_4}\left(\frac{u_3}{\rho_n}\right)e^{i\frac{u_3\lambda_2}{\rho_n}}
1_{\left((u_1-1)b_n,u_1 b_n)\right]}(x_1)
1_{\left(\frac{(u_2-1)}{\rho_n},\frac{u_2}{\rho_n}\right]}(x_2)1_{\left(\frac{(u_3-1)}{\rho_n},\frac{u_3}{\rho_n}\right]}(x_3).
\end{split}
\end{equation*}
So that
$$nb_nT_{11}(\lambda_1,\lambda_2)=\frac{1}{(4\pi)^{2}}\int_{-\infty}^{\infty}\int_{-\infty}^{\infty}\int_{-\infty}^{\infty}S_n(\lambda_1,\lambda_2,x_1,x_2,x_3)dx_1dx_2dx_3.$$

Define $u_{1n}(x_1)$, $u_{2n}(x_2)$ and $u_{3n}(x_3)$ as the
smallest integers greater than or equal to $x_1/b_n$, $\rho_n x_2$
and $\rho_nx_3$, respectively. Thus, $(x_1,x_2,x_3) \in
(b_nu_{1n-1}(x_1),b_n
u_{1n}(x_1)]\times\left(\frac{u_{2n-1}(x_2)}{\rho_n},\frac{u_{2n}(x_2)}{\rho_n}\right]\times
\left(\frac{u_{3n-1}(x_3)}{\rho_n},\frac{u_{3n}(x_3)}{\rho_n}\right]$
and $b_n u_{1n}(x_1)\!\rightarrow\! x_1
,\frac{u_{2n}(x_2)}{\rho_n}\!\rightarrow\! x_2,
\frac{u_{3n}(x_3)}{\rho_n}\! \rightarrow\! x_3$ as $n\rightarrow
\infty$. Since $nb_n\rightarrow\infty$ and $b_n\rho_n\rightarrow0$
as $n\rightarrow\infty$, we have, for any point $(x_1,x_2,x_3)\in
\mathbb{R}^3$ and large enough $n$, the inequalities
$-\frac{nb_n-x_1}{b_n\rho_n}<x_3<\frac{nb_n-x_1}{b_n\rho_n}$, i.e.,
$-n+1-u_{1n}(x_1)< u_{3n}(x_3)< n-1-u_{1n}(x_1)$. Thus, for
sufficiently large $n$, we have
\begin{equation}
\begin{split}
&\hskip-15pt
S_n(\lambda_1,\lambda_2,x_1,x_2,x_3)\\&=U_n(u_{1n}(x_1),u_{2n}(x_2),u_{3n}(x_3))
K(b_n u_{1n}(x_1))K(b_n(u_{1n}(x_1)-u_{2n}(x_2)+u_{3n}(x_3)))
\\&~~\times C_{a_1a_3}\left(\frac{u_{2n}(x_2)}{\rho_n}\right)e^{-i\frac{u_{2n}(x_2)\lambda_2}{\rho_n}} C_{a_2a_4}\left(\frac{u_{3n}(x_3)}{\rho_n}\right)e^{i\frac{u_{3n}(x_3)\lambda_2}{\rho_n}}.
\end{split}
\label{largeq}
\end{equation}

Observe that, under Assumptions 1,3 and~4, the function
$S_n(\lambda_1,\lambda_2,x_1,x_2,x_3)$ converges to the function
$S(\cdot)$, defined over $E'\times\mathbb{R}^3$ by
\begin{equation*}
S(\lambda_1,\lambda_2,x_1,x_2,x_3)=K^{2}(x_1)C_{a_1a_3}\left(x_2\right)e^{-ix_2\lambda_2}
C_{a_2a_4}\left(x_3\right)e^{ix_3\lambda_2}.
\end{equation*}
Observe also that
$\int_{-\infty}^{\infty}\int_{-\infty}^{\infty}\int_{-\infty}^{\infty}S(\lambda_1,\lambda_2,x_1,x_2,x_3)dx_1dx_2dx_3$
is a continuous function in $(\lambda_1,\lambda_2)$. As in the proof
of Theorem~1, we prove the convergence of the left hand side of
\eqref{largeq} uniformly on $E'$, by showing that for any sequence
$(\lambda_{1n},\lambda_{2n}) \rightarrow (\lambda_{1},\lambda_{2})$,
\begin{equation*}
\begin{split}
&\hskip-50pt\lim_{n\rightarrow\infty}\int_{-\infty}^{\infty}\int_{-\infty}^{\infty}\int_{-\infty}^{\infty}
S_n(\lambda_{1n},\lambda_{2n},x_1,x_2,x_3)dx_1dx_2dx_3\\
&=\int_{-\infty}^{\infty}\int_{-\infty}^{\infty}\int_{-\infty}^{\infty}S(\lambda_{1},\lambda_{2},x_1,x_2,x_3)dx_1dx_2dx_3.\\
\end{split}
\end{equation*}
for $(\lambda_{1n},\lambda_{2n}),(\lambda_{1},\lambda_{2}) \in E' $.
The latter convergence follows, through Assumption~1 and~2 and the
DCT, from the dominance
$$|S_n(\lambda_{1n},\lambda_{2n},x_1,x_2,x_3)|\le M K_1(x_1)g_{a_1a_3}\left(x_2\right)g_{a_2a_4}\left(x_3\right).$$
and the convergence of the integrand, which holds because of the
continuity of $C_{a_1a_3}(\cdot)$, $C_{a_2a_4}(\cdot)$ and the
kernel and the exponential functions. Hence, $nb_nT_{11}(\cdot)$
converges as stated uniformly on the compact set~$E'.\hfill\Box$

\bigskip{\sc Lemma~2.} {\em For $\lambda_1-\lambda_2\ne0$, the function
$nb_nT_{11}(\lambda_1,\lambda_2)$ converges to zero. The convergence
is uniform on any  compact subset of the set $E_1$ given by
$$E=\{(\lambda_1,\lambda_2):\lambda_1-\lambda_2\ne0,~-\infty<\lambda_1,\lambda_2<\infty\}.$$}

{\sc Proof of Lemma~2.} Let $E'$ be any compact subset of the set $E$. Consider the simple function
$S_n(\cdot)$, defined over $E'\times\mathbb{R}^3$ by
\begin{equation*}
\begin{split}
\hskip-10pt&S_n(\lambda_1,\lambda_2,x_1,x_2,x_3)\\&=\sum_{u_1=-(n-1)}^{(n-1)}\sum_{u_2=-(n-1)}^{(n-1)}\sum_{u_3=-(n-1)-u_1}^{(n-1)-u_1}U_n(u_1,u_2,u_3)
K(b_n u_1)K(b_n(u_1-u_2+u_3))e^{-i\frac{u_1
(\lambda_1-\lambda_2)}{\rho_n}}
\\&~~\times  C_{a_1a_3}\left(\frac{u_2}{\rho_n}\right)e^{-i\frac{u_2\lambda_2}{\rho_n}} C_{a_2a_4}\left(\frac{u_3}{\rho_n}\right)e^{i\frac{u_3\lambda_2}{\rho_n}}
1_{\left((u_1-1)b_n,u_1 b_n)\right]}(x_1)
1_{\left(\frac{(u_2-1)}{\rho_n},\frac{u_2}{\rho_n}\right]}(x_2)1_{\left(\frac{(u_3-1)}{\rho_n},\frac{u_3}{\rho_n}\right]}(x_3).
\end{split}
\end{equation*}
So that
$$nb_nT_{11}(\lambda_1,\lambda_2)=\frac{1}{(4\pi)^{2}}\int_{-\infty}^{\infty}\int_{-\infty}^{\infty}\int_{-\infty}^{\infty}S_n(\lambda_1,\lambda_2,x_1,x_2,x_3)dx_1dx_2dx_3.$$
An argument similar to that used in the proof of Lemma 1 shows that
for $(x_1,x_2,x_3)\in \mathbb{R}^3$ and sufficiently large $n$,
\begin{equation*}
\begin{split}
&\hskip-15pt
S_n(\lambda_1,\lambda_2,x_1,x_2,x_3)\\&=U_n(u_{1n}(x_1),u_{2n}(x_2),u_{3n}(x_3))
K(b_n u_{1n}(x_1))K(b_n(u_{1n}(x_1)-u_{2n}(x_2)+u_{3n}(x_3)))
\\&~~\times e^{-i\frac{u_{1n}(x_1)(\lambda_1-\lambda_2)}{\rho_n}}C_{a_1a_3}\left(\frac{u_{2n}(x_2)}{\rho_n}\right)e^{-i\frac{u_{2n}(x_2)\lambda_2}{\rho_n}} C_{a_2a_4}\left(\frac{u_{3n}(x_3)}{\rho_n}\right)e^{i\frac{u_{3n}(x_3)\lambda_2}{\rho_n}}.
\end{split}
\end{equation*}
where $u_{1n}(x_1)$, $u_{2n}(x_2)$ and $u_{3n}(x_3)$ are the
smallest integers greater than or equal to $x_1/b_n$, $\rho_nx_2$
and $\rho_nx_3$, respectively.

For obtaining the uniform convergence of
$nb_nT_{11}(\lambda_1,\lambda_2)$, consider
\begin{align}\label{twostep}
&\hskip-20pt\sup_{(\lambda_1,\lambda_2) \in
E'}\left|\int_{-\infty}^{\infty}\int_{-\infty}^{\infty}\int_{-\infty}^{\infty}S_n(\lambda_1,\lambda_2,x_1,x_2,x_3)dx_1dx_2dx_3\right|\notag\\
\le&\sup_{(\lambda_1,\lambda_2) \in
E'}\int_{-\infty}^{\infty}\int_{-\infty}^{\infty}\int_{-\infty}^{\infty}|S_n(\lambda_1,\lambda_2,x_1,x_2,x_3)-g_n(\lambda_1,\lambda_2,x_1,x_2,x_3)|dx_1dx_2dx_3\notag\\
+&\sup_{(\lambda_1,\lambda_2) \in
E'}\left|\int_{-\infty}^{\infty}\int_{-\infty}^{\infty}\int_{-\infty}^{\infty}g_n(\lambda_1,\lambda_2,x_1,x_2,x_3)dx_1dx_2dx_3\right|,
\end{align}
where the function $g_n(\cdot)$ is defined over
$E'\times\mathbb{R}^3$ by
$$g_n(\lambda_1,\lambda_2,x_1,x_2,x_3)
=K^{2}(x_1)e^{-i\frac{x_1 (\lambda_1-\lambda_2)}{b_n\rho_n}} C_{a_1a_3}\left(x_2\right)e^{-ix_2\lambda_2} C_{a_2a_4}\left(x_3\right)e^{ix_3\lambda_2}.$$ %
We will show the uniform convergence of the right hand side of
\eqref{twostep} by considering the two terms separately. For the
first term, we follow the route taken in the proof of Theorem~1,
i.e., show that for any sequence
$(\lambda_{1n},\lambda_{2n})\rightarrow(\lambda_1,\lambda_2)$,
$$\lim_{n\rightarrow\infty}\int_{-\infty}^{\infty}\int_{-\infty}^{\infty}\int_{-\infty}^{\infty}
|S_n(\lambda_{1n},\lambda_{2n},x_1,x_2,x_3)-g_n(\lambda_{1n},\lambda_{2n},x_1,x_2,x_3)|dx_1dx_2dx_3=0$$
for $(\lambda_{1n},\lambda_{2n}),(\lambda_1,\lambda_2)\in E'$. For
this purpose, we write the above integral as
\begin{align}\label{2part}
&\hskip-15pt\int_{-\infty}^{\infty}\int_{-\infty}^{\infty}\int_{-\infty}^{\infty}
|S_n(\lambda_{1n},\lambda_{2n},x_1,x_2,x_3)-g_n(\lambda_{1n},\lambda_{2n},x_1,x_2,x_3)|dx_1dx_2dx_3\notag\\
\le&\int_{-\infty}^{\infty}\int_{-\infty}^{\infty}\int_{-\infty}^{\infty}
|S_n(\lambda_{1n},\lambda_{2n},x_1,x_2,x_3)-G_n(\lambda_{1n},\lambda_{2n},x_1,x_2,x_3)|dx_1dx_2dx_3\notag\\
&+\int_{-\infty}^{\infty}\int_{-\infty}^{\infty}\int_{-\infty}^{\infty}
|G_n(\lambda_{1n},\lambda_{2n},x_1,x_2,x_3)-g_n(\lambda_{1n},\lambda_{2n},x_1,x_2,x_3)|dx_1dx_2dx_3,
\end{align}
where the function $G_n(\cdot)$ is defined over
$E'\times\mathbb{R}^3$ by
$$G_n(\lambda_1,\lambda_2,x_1,x_2,x_3)
=K^{2}(x_1)e^{-i\frac{u_{1n}(x_1)b_n
(\lambda_1-\lambda_2)}{b_n\rho_n}}
C_{a_1a_3}\left(x_2\right)e^{-ix_2\lambda_2}
C_{a_2a_4}\left(x_3\right)e^{ix_3\lambda_2}.$$

Now observe that
\begin{align*}
&\hskip-50pt
|S_n(\lambda_{1n},\lambda_{2n},x_1,x_2,x_3)-G_n(\lambda_{1n},\lambda_{2n},x_1,x_2,x_3)|\\
&\le M \left|e^{-i\frac{u_{1n}(x_1)b_n
(\lambda_{1n}-\lambda_{2n})}{b_n\rho_n}}\alpha_n(\lambda_{1n},\lambda_{2n},x_1,x_2,x_3)\right|,
\end{align*}
where
\begin{align*}
&\hskip-20pt \alpha_n(\lambda_{1n},\lambda_{2n},x_1,x_2,x_3)\\
&= U_n(u_{1n}(x_1),u_{2n}(x_2),u_{3n}(x_3))
K(b_n u_{1n}(x_1))K(b_n(u_{1n}(x_1)-u_{2n}(x_2)+u_{3n}(x_3)))\\
&~~~\times C_{a_1a_3}\left(\frac{u_{2n}(x_2)}{\rho_n}\right)e^{-i\frac{u_{2n}(x_2)\lambda_{2n}}{\rho_n}} %
C_{a_2a_4}\left(\frac{u_{3n}(x_3)}{\rho_n}\right)e^{i\frac{u_{3n}(x_3)\lambda_{2n}}{\rho_n}}\\
&~~~-K^{2}(x_1)C_{a_1a_3}\left(x_2\right)e^{-ix_2\lambda_{2n}}C_{a_2a_4}\left(x_3\right)e^{ix_3\lambda_{2n}}.
\end{align*}
Since $\alpha_n(\lambda_n,x,t,t^{'})\rightarrow 0$ as $n\rightarrow
\infty$, we have
$$\lim_{n\rightarrow\infty}|S_n(\lambda_{1n},\lambda_{2n},x_1,x_2,x_3)-G_n(\lambda_{1n},\lambda_{2n},x_1,x_2,x_3)|=0$$
Since from Assumption 1 and 2, we have the dominance
\begin{align*}
|S_n(\lambda_{1n},\lambda_{2n},x_1,x_2,x_3)-G_n(\lambda_{1n},\lambda_{2n},x_1,x_2,x_3)|\le
2M K_1(x_1)g_{a_1a_3}(x_2)g_{a_2a_4}(x_2).
\end{align*}
By applying the DCT, we have
$$\lim_{n\rightarrow\infty}\int_{-\infty}^{\infty}\int_{-\infty}^{\infty}\int_{-\infty}^{\infty}|S_n(\lambda_{1n},\lambda_{2n},x_1,x_2,x_3)-G_n(\lambda_{1n},\lambda_{2n},x_1,x_2,x_3)|dx_1dx_2dx_3=0.$$

Turning to the second term on the right hand side of (\ref{2part}),
observe that for any fixed $x_1$,
$$\left|e^{-i\frac{u_{1n}(x_1)b_n (\lambda_{1n}-\lambda_{2n})}{b_n\rho_n}}-e^{-i\frac{x_1 (\lambda_{1n}-\lambda_{2n})}{b_n\rho_n}}\right|\le \frac{\lambda_{1n}-\lambda_{2n}}{\rho_n}.$$
Thus,
\begin{align*}
\hskip-30pt&|G_n(\lambda_{1n},\lambda_{2n},x_1,x_2,x_3)-g_n(\lambda_{1n},\lambda_{2n},x_1,x_2,x_3)|
\\&\leq M^{2}g_{a_1a_3}(0)g_{a_2a_4}(0)\left|e^{-i\frac{u_{1n}(x_1)b_n (\lambda_{1n}-\lambda_{2n})}{b_n\rho_n}}-e^{-i\frac{x_1 (\lambda_{1n}-\lambda_{2n})}{b_n\rho_n}}\right|
\le
M^{2}g_{a_1a_3}(0)g_{a_2a_4}(0)\frac{\lambda_{1n}-\lambda_{2n}}{\rho_n},
\end{align*}
and so
$$\lim_{n\rightarrow\infty}|G_n(\lambda_{1n},\lambda_{2n},x_1,x_2,x_3)-g_n(\lambda_{1n},\lambda_{2n},x_1,x_2,x_3)|= 0.$$
From Assumption 1 and 2, we have the dominance
\begin{align*}
|G_n(\lambda_{1n},\lambda_{2n},x_1,x_2,x_3)-g_n(\lambda_{1n},\lambda_{2n},x_1,x_2,x_3)|\le
2M K_1(x_1)g_{a_1a_3}(x_2)g_{a_2a_4}(x_2).
\end{align*}
which leads us, through another use of the DCT, to the convergence
of the second integral of (\ref{2part}). This establishes that the
first term on the right hand side of (\ref{twostep}) converges to~0.
We only have to deal with the second term. Let
$$s_n(\lambda_1,\lambda_2)=\int_{-\infty}^{\infty}\int_{-\infty}^{\infty}\int_{-\infty}^{\infty}g_n(\lambda_{1},\lambda_{2},x_1,x_2,x_3)dx_1dx_2dx_3.$$
In order to establish the uniform convergence of $s_n(\cdot)$ over
$E'$, it is enough to show that
$s_n(\lambda_{1n},\lambda_{2n})\rightarrow 0$ for any sequence
$(\lambda_{1n},\lambda_{2n})\rightarrow(\lambda_{1},\lambda_{2})$,
where $(\lambda_{1n},\lambda_{2n}),(\lambda_{1},\lambda_{2})\in E'$.
By using the Reimann-Lebesgue lemma, we have
$s_n(\lambda_1,\lambda_2)\rightarrow 0$. Thus, the second term on
the right hand side of (\ref{twostep}) also converges to~0. Hence,
$nb_nT_{11}(\lambda_1,\lambda_2)$ converges to~0 uniformly on $E'$
as $n\rightarrow\infty.\hfill\Box$

\medskip{\sc Proof of Theorem~3.} $cum(\widehat{\phi}_{a_1a_2}(\lambda_1),
\widehat{\phi}_{a_3a_4}(\lambda_2),\ldots,\widehat{\phi}_{a_{2L-1}a_{2L}}(\lambda_L))$
can be written as
\begin{equation}\label{jointcum}
\begin{split}
\hskip-15pt&cum(\widehat{\phi}_{a_1a_2}(\lambda_1),\widehat{\phi}_{a_3a_4}(\lambda_2),\ldots,\widehat{\phi}_{a_{2L-1}a_{2L}}(\lambda_L))\\
&=\frac{1}{(\pi n
\rho_n)^{L}}\sum_{t_1=1}^{n}\sum_{t_2=1}^{n}\ldots\sum_{t_{2L-1}=1}^{n}\sum_{t_{2L}=1}^{n}K(b_n(t_1-t_2))
\ldots K(b_n(t_{2L-1}-t_{2L}))e^{-\frac{i(t_1-t_2)\lambda_1}{\rho_n}}\times\cdots\\
&~~\times
e^{-\frac{i(t_{2L-1}-t_{2L})\lambda_L}{\rho_n}}cum\left(X_{a_1}\left(\frac{t_1}{\rho_n}\right)X_{a_2}
\left(\frac{t_2}{\rho_n}\right),\ldots,
X_{a_{2L-1}}\left(\frac{t_{2L-1}}{\rho_n}\right)X_{a_{2L}}\left(\frac{t_{2L}}{\rho_n}\right)\right)
\end{split}
\end{equation}
It follows that
\begin{equation*}
\begin{split}
&\hskip-30pt|cum(\widehat{\phi}_{a_1a_2}(\lambda_1),\widehat{\phi}_{a_3a_4}(\lambda_2),\ldots,\widehat{\phi}_{a_{2L-1}a_{2L}}(\lambda_L))|\\&
\le\frac{1}{(n
\rho_n)^{L}}\sum_{t_1=1}^{n}\sum_{t_2=1}^{n}\ldots\sum_{t_{2L-1}=1}^{n}\sum_{t_{2L}=1}^{n}|K(b_n
(t_1-t_2))\ldots K(b_n (t_{2L-1}-t_{2L}))|\\
&\qquad\times\left|cum\left(X_{a_1}\left(\frac{t_1}{\rho_n}\right)X_{a_2}\left(\frac{t_2}{\rho_n}\right),\ldots,
X_{a_{2L-1}}\left(\frac{t_{2L-1}}{\rho_n}\right)X_{a_{2L}}\left(\frac{t_{2L}}{\rho_n}\right)\right)\right|
\end{split}
\end{equation*}
Now
\begin{equation*}
\begin{split}
&\hskip-30pt
cum\left(X_{a_1}\left(\frac{t_1}{\rho_n}\right)X_{a_2}\left(\frac{t_2}{\rho_n}\right),\ldots,
X_{a_{2L-1}}\left(\frac{t_{2L-1}}{\rho_n}\right)X_{a_{2L}}\left(\frac{t_{2L}}{\rho_n}\right)\right)\\
&=\sum_{\bn}C_{a_{j_{11}}a_{j_{12}}\ldots
a_{j_{1k_1}}}\left(\frac{t_{j_{11}}-t_1'}{\rho_n},\ldots,\frac{t_{j_{1,k_1-1}}-t_1'}{\rho_n} \right)\cdots\\
& \qquad\times C_{a_{j_{P1}}a_{j_{P2}}\ldots a_{j_{Pk_P}}}
\left(\frac{t_{j_{P1}}-t_P'}{\rho_n},\ldots,\frac{t_{j_{P,k_P-1}}-t_P'}{\rho_n} \right)%
\end{split}
\end{equation*}
where the summation is over all {\it indecomposable} (Brillinger,
2001; Leonov and Shiryayev, 1959) partitions
$\bn=(\nu_1,\ldots,\nu_P)$, such that
$\nu_p=(j_{p1},\ldots,j_{pk_p})$, $p=1,\ldots,P$, of the table
$$
\begin{tabular}{cc}
1&2\\
3&4\\
\vdots &\vdots\\
2L-1&2L
\end{tabular}
$$
and $t_p'=t_{j_{pk_p}}$, $p=1,\ldots,P$. Since the partition $\bn$
is indecomposable, we have
$$t_{j_{pl}}-t_p'\ne t_{2m}-t_{2m-1};\ l=1,\ldots,k_p;\ p=1,\ldots,P;\ m=1,\ldots,L.$$
Define
$$u_{j_{pl}}=t_{j_{pl}}-t_p';\ l=1,\ldots,k_p; p=1,\ldots,P.$$
Note that $u_{j_{pk_p}}=0$ for $p=1,\ldots,P$. Then the joint
cumulant of $(\widehat{\phi}_{a_1a_2}(\lambda_1)$,
$\widehat{\phi}_{a_3a_4}(\lambda_2)$, $\ldots$,
$\widehat{\phi}_{a_{2L-1}a_{2L}}(\lambda_L))$ given by
\eqref{jointcum} is absolutely bounded by
\begin{equation}\label{2ndbound}
\begin{split}
&\frac{1}{(n \rho_n)^{L}}\sum_{\bn} %
\sum_{t_1'=1}^{n}\sum_{u_{j_{11}}=-(t_1'-1)}^{n-t_1'}\ldots\sum_{u_{j_{1,k_1-1}}=-(t_1'-1)}^{n-t_1'}\ldots %
\sum_{t_P'=1}^{n}\sum_{u_{j_{P1}}=-(t_P'-1)}^{n-t_P'}\ldots\sum_{u_{j_{P,k_P-1}}=-(t_P'-1)}^{n-t_P'}\\
&~~~\left|K[b_n (u_1    +t_{p_1}'     -u_2-t_{p_2}')]\right|\cdots
    \left|K[b_n(u_{2L-1}+t_{p_{2L-1}}'-u_{2L}-t_{p_{2L}}')]\right|\\
&\qquad\times\left|
C_{a_{j_{11}}a_{j_{12}}\ldots a_{j_{1k_1}}}\left(\frac{u_{j_{11}}}{\rho_n},\ldots,\frac{u_{j_{1,k_1-1}}}{\rho_n} \right)%
\ldots %
C_{a_{j_{P1}}a_{j_{P2}}\ldots a_{j_{Pk_P}}}\left(\frac{u_{j_{P1}}}{\rho_n},\ldots,\frac{u_{j_{P,k_P-1}}}{\rho_n} \right)\right|,%
\end{split}
\end{equation}
where $p_m$ is that member of the set $\{1,2,\ldots,P\}$ which
satisfies $t_m\in\nu_{p_m}$ for $m=1,\ldots,L$.

We will now show that the set
$A=\{t_{p_1}'-t_{p_2}',\ldots,t_{p_{2L-1}}'-t_{p_{2L}}'\}$ has $P-1$
linearly independent elements. Note that the set $A$ consists of
differences of pairs of elements of the set
$\{t_1',t_2',\ldots,t_P'\}$. So the set $A$ can have at most $P-1$
linearly independent differences. Suppose that the set $A$ has
exactly $P-j$ linearly independent differences for some $j\ge1$.
Denote the $P-j$ independent differences of the set $A$ by
$$A_1=\left\{t_{p_{2k_1-1}}'-t_{p_{2k_1}}',t_{p_{2k_2-1}}'-t_{p_{2k_2}}',
\ldots,t_{p_{2k_{P-j}-1}}'-t_{p_{2k_{P-j}}}'\right\},$$ %
where $k_1,\ldots,k_{P-j}\in\{1,2,\ldots,L\}$. Let, if possible,
$j>1$, and consider a difference $t_{l_1}'-t_{l_2}'$ for
$l_1,l_2\in\{1,2,\ldots,P\}$ which is linearly independent of the
elements of the set $A_1$. Since the partition $\bn$ is
indecomposable, the sets $\nu_{l_1}$ and $\nu_{l_2}$ {\it
communicate} (Leonov and Shiryayev, 1959). Therefore, there exists
an index set $\{s_1,s_2,\ldots,s_r\}$ with $r\ge2$, which is a
proper subset of $\{1,2,\ldots,P\}$, such that $s_1=l_1$, $s_r=l_2$
and the pairs
$(\nu_{s_1},\nu_{s_2}),(\nu_{s_2},\nu_{s_3}),\ldots,(\nu_{s_{r-1}},\nu_{s_r})$
are {\it hook} (Leonov and Shiryayev, 1959). Consequently, there
exist indices $j_1,\ldots,j_{r-1}\in\{1,\ldots,L\}$ such that for
$m=1,\ldots,r-1$, one of the points $t_{2j_m-1}$ and $t_{2j_m}$
belongs to $\nu_{s_m}$ and the other belongs to $\nu_{s_{m+1}}$. It
follows that for $m=1,\ldots,r-1$, $(t_{p_{2j_m-1}}'-t_{p_{2j_m}}')$
is in $A$, and hence, they can be written as linear combinations of
the members of $A_1$. Note that for $m=1,\ldots,r-1$,
$(t_{s_{m-1}}'-t_{s_m}')$ is equal to either
$(t_{p_{2j_m-1}}'-t_{p_{2j_m}}')$ or
$-(t_{p_{2j_m-1}}'-t_{p_{2j_m}}')$. Thus,
$$t_{l_1}'-t_{l_2}'=t_{s_1}'-t_{s_r}'
=(t_{s_1}'-t_{s_2}')+(t_{s_2}'-t_{s_3}')+\cdots+(t_{s_{r-1}}'-t_{s_r}')$$ %
can be written as a linear combination of the members of $A_1$. This
fact contradicts the assumption that $t_{l_1}'-t_{l_2}'$ is linearly
independent of the elements of the set $A_1$. Therefore, $j$ cannot
be larger than 1. This proves that the set $A$ cannot contain fewer
than $P-1$ linearly independent differences.

Consider the $P-1$ linearly independent elements of the set $A_1$,
where $j=1$, and define
\begin{eqnarray*}
v_1&=&u_{2k_1-1} + t_{p_{2k_1-1}}' - u_{2k_1} - t_{p_{2k_1}}',\\
&\vdots&\\
v_{P-1}&=&u_{2k_{P-1}-1} + t_{p_{2k_{P-1}-1}}' - u_{2k_{P-1}} -
t_{p_{2k_{P-1}}}'.%
\end{eqnarray*}
Using the above transformation, and by replacing the $P$ sums over
indices $t_1',\ldots,t_P'$ by $P-1$ sums overs the indices
$v_1,\ldots,v_{P-1}$, we find that the joint cumulant given in
(\ref{2ndbound}) is bounded from above by
\begin{equation}\label{3rdbound}
\begin{split}
&\frac1{n^{L-1}\rho_n^L}\sum_{\bn}M^{L-P+1}\sum_{u_{j_{11}}=-(n-1)}^{n-1}\ldots\sum_{u_{j_{1,k_1-1}}
=-(n-1)}^{n-1}\ldots \sum_{u_{j_{P1}}=-(n-1)}^{n-1}\ldots\sum_{u_{j_{P,k_{P}-1}}=-(n-1)}^{n-1}\\
&\sum_{v_1=-3n}^{3n}\ldots\sum_{v_{P-1}=-3n}^{3n}\left|K(b_nv_1)\right|\ldots\left|K(b_nv_{P-1})\right|\\
&\times\left|C_{a_{j_{11}}a_{j_{12}}\ldots
a_{j_{1k_1}}}\left(\frac{u_{j_{11}}}{\rho_n},\ldots,
\frac{u_{j_{1,k_1-1}}}{\rho_n} \right)\right|\ldots
\left|C_{a_{j_{P1}}a_{j_{P2}}\ldots
a_{j_{Pk_P}}}\left(\frac{u_{j_{P1}}}{\rho_n},\ldots,
\frac{u_{j_{P,k_P-1}}}{\rho_n} \right)\right|.
\end{split}
\end{equation}
The above simplification has been made by taking into account the
upper bound for $L-P+1$ copies of $K(\cdot)$ and conservative
estimates of the ranges of summation of $v_1,\ldots,v_{P-1}$. Now
one can rewrite the expression in (\ref{3rdbound}) as follows.
\begin{equation}\label{4thbound}
\begin{split}
&\sum_{\bn}M^{L-P+1}~\frac{(\rho_nb_n)^{L-P}}{(nb_n)^{L-1}}\left[\sum_{v_1=-3n}^{3n}K(b_n v_1)b_n\right]\ldots%
\left[\sum_{v_{P-1}=-3n}^{3n}K(b_n v_{P-1})b_n\right]\\
&\times\left\{\sum_{u_{j_{11}}=-(n-1)}^{n-1}\ldots\sum_{u_{j_{1,k_1-1}}=-(n-1)}^{n-1}%
\left|C_{a_{j_{11}}a_{j_{12}}\ldots a_{j_{1k_1}}}%
\left(\frac{u_{j_{11}}}{\rho_n},\ldots, \frac{u_{j_{1,k_1-1}}}{\rho_n} \right)\right|%
\frac{1}{\rho_n^{k_1-1}}\right\}\ldots\\
&\times\left\{\sum_{u_{j_{P1}}=-(n-1)}^{n-1}\ldots\sum_{u_{j_{P,k_{P}-1}}=-(n-1)}^{n-1}%
\left|C_{a_{j_{P1}}a_{j_{P2}}\ldots
a_{j_{Pk_P}}}\left(\frac{u_{j_{P1}}}{\rho_n},\ldots,
\frac{u_{j_{P,k_P-1}}}{\rho_n}
\right)\right|\frac{1}{\rho_n^{k_P-1}}\right\}.
\end{split}
\end{equation}

Consider the simple function $S_n(\cdot)$ defined over $\mathbb{R}$
by
$$S_n(x)=\sum_{v_1=-3n}^{3n}K(b_n v_1)1_{\left(b_nv_1-1,b_nv_1\right]}(x).$$
Note that $\int_{-\infty}^{\infty}S_n(x)dx=\sum_{v_1=-3n}^{3n}K(b_n
v_1)b_n$, and from Assumption 2 we have the dominance $S_{n}(x)\le
K_1(x)$. By applying the DCT, we have
$$\sum_{v_1=-3n}^{3n}K(b_n v_1)b_n\rightarrow\int_{-\infty}^{\infty}|K(x)|dx.$$
This fact establishes the convergence of the sums over
$v_1,\ldots,v_{P-1}$.

Consider the simple function $T_n(\cdot)$ defined over
$\mathbb{R}^{k_1-1}$ by
\begin{equation*}
\begin{split}
T_n(x_1,x_2,\ldots,x_{k_1-1})=&\sum_{u_{j_{11}}=-(n-1)}^{n-1} \ldots
\sum_{u_{j_{1,k_{1}-1}}=-(n-1)}^{n-1}C_{a_{j_{11}}a_{j_{12}}\ldots
a_{j_{1k_1}}}\left(\frac{u_{j_{11}}}{\rho_n},\ldots,
\frac{u_{j_{1,k_1-1}}}{\rho_n}\right)\\& \times 1_{\left(
\frac{u_{j_{11}}-1}{\rho_n},\frac{u_{j_{11}}}{\rho_n}\right]}(x_1)\ldots
1_{\left(\frac{u_{j_{1,k_{1}-1}}-1}{\rho_n},\frac{u_{j_{1,k_{1}-1}}}{\rho_n}\right]}(x_{k_1-1}).
\end{split}
\end{equation*}
Note that
\begin{equation*}
\begin{split}
&\int_{-\infty}^{\infty}\ldots\int_{-\infty}^{\infty}T_n(x_1,\ldots,x_{k_1-1})dx_1\ldots
dx_{k_1-1}\\&
\hskip10pt=\sum_{u_{j_{11}}=-(n-1)}^{n-1}\ldots\sum_{u_{j_{1,k_{1}-1}}=-(n-1)}^{n-1}
\left|C_{a_{j_{11}}a_{j_{12}}\ldots a_{j_{1k_1}}}
\left(\frac{u_{j_{11}}}{\rho_n},\ldots,
\frac{u_{j_{1,k_1-1}}}{\rho_n}
\right)\right|\frac{1}{\rho_n^{k_1-1}}.
\end{split}
\end{equation*}
From Assumption~5A, we have that the function $T_n(\cdot)$ is
bounded by an integrable function. Thus, by applying the DCT, we
have
\begin{equation*}
\begin{split}
&\lim_{n\rightarrow\infty}
\sum_{u_{j_{11}}=-(n-1)}^{n-1}\ldots\sum_{u_{j_{1,k_{1}-1}}=-(n-1)}^{n-1}
\left|C_{a_{j_{11}}a_{j_{12}}\ldots a_{j_{1k_1}}}
\left(\frac{u_{j_{11}}}{\rho_n},\ldots,
\frac{u_{j_{1,k_1-1}}}{\rho_n}
\right)\right|\frac{1}{\rho_n^{k_1-1}}\\
&\hskip20pt=\int_{-\infty}^{\infty}\ldots\int_{-\infty}^{\infty}
\left|C_{a_{j_{11}}a_{j_{12}}\ldots
a_{j_{1k_1}}}\left(x_1,\ldots,x_{k_1-1}\right)\right|dx_1\ldots
dx_{k_1-1}.
\end{split}
\end{equation*}
Likewise, we have the convergence for the remaining $P-1$ sets of
summations. Using these above convergence results, the upper bound
of (\ref{2ndbound}) given in (\ref{4thbound}) can be written as
$$\sum_{\bn}\frac{(\rho_nb_n)^{L-P}}{(nb_n)^{L-1}} d_{\bn},$$
where $d_{\bn}$ are appropriate constants. The summation is over the
finite number of indecomposable partitions, and the worst-case value
of the partition size $P$ is $L$. Therefore, the upper bound is
$O\left((nb_n)^{-(L-1)}\right)$. This Completes the proof of
Theorem~3.\hfill$\Box$

\medskip{\sc Proof of Theorem~4.} Note that the first moment of the random
vector on the left hand side of (\ref{meancorrected}) is zero and
the second moment converges in accordance with Theorem~2. Further,
$$cum(c_1(Y_1-d_1),c_1(Y_2-d_2),\ldots,c_J(Y_J-d_J))=c_1c_2\cdots c_J\times cum(Y_1,Y_2,\ldots,Y_J),$$
for any set of constants $c_1,\ldots,c_J,d_1,\ldots,d_J$. From the
above fact and Theorem~3, for all $k>2$, the absolute value of the
$k$th order joint cumulant of the random vector on the left hand
side of (\ref{meancorrected}) is bounded from above by an
$O((nb_n)^{k/2-k+1})$ term. According to Assumption~3, this upper
bound tends to 0 as $n$ tends to infinity. This completes the
proof.\hfill$\Box$

\medskip{\sc Proof of Theorem~5.} The result can be proved along the lines
of the proof of Theorem~3 of Srivastava and Sengupta
(2010).\hfill$\Box$

\medskip{\sc Proof of Theorem~6.} The weak convergence of the first term on
the right hand side of (\ref{asym1}) follows from Theorem~4. On the
other hand, the second term can be written, in view of Theorem~5, as
\begin{equation}\label{asymrate}
\sqrt{nb_n}\left(E[\widehat{\phi}_{a_1a_2}(\lambda_1)]-\phi_{a_1a_2}(\lambda_1)\right)
=\sqrt{nb_n}\left(O\left((\rho_nb_n)^{q}\right)+O\left(\frac{\rho_n}{n}\right)+O\left(\frac{1}{\rho_n^{p}}\right)\right).
\end{equation}
Under Assumption 3, %
$$\lim_{n\rightarrow
\infty}\sqrt{nb_n}\rho_n^qb_n^q=0\Rightarrow
\lim_{n\rightarrow\infty}\sqrt{nb_n}\frac{\rho_n}{n}=0.$$ %
Therefore, under Assumptions 3 and 4A, the right hand side of
(\ref{asymrate}) goes to zero as $n\rightarrow\infty$. This
completes the proof.\hfill$\Box$

\medskip{\sc Proof of Theorem~7.} Note that under Assumption~4A, we have
\begin{align}\label{2nd}
\lim_{n\rightarrow \infty}\sqrt{nb_n}\frac{1}{\rho_n^p}=0
\end{align}
and
\begin{eqnarray}
&\displaystyle\lim_{n\rightarrow \infty}\sqrt{nb_n}\rho_n^qb_n^q=0%
\quad&\Leftrightarrow \lim_{n\rightarrow\infty}\left(nb_n\right)^{\frac1{2q}}b_n\rho_n=0\nonumber\\
\Leftrightarrow &\displaystyle\lim_{n\rightarrow\infty}\left(nb_n\right)^{1+\frac1{2q}}\frac{\rho_n}n=0%
&\Leftrightarrow \lim_{n\rightarrow \infty}\sqrt{nb_n}
\left(\frac{\rho_n}{n}\right)^{\frac{q}{1+2q}}=0. \label{1st}
\end{eqnarray}

From (\ref{2nd}) and (\ref{1st}), we have
\begin{eqnarray}
\frac{1}{\sqrt{nb_n}}&=&o\left(\left(\frac{\rho_n}{n}\right)^{\frac{q}{1+2q}}\right),\label{rate1}\\
\mbox{and
}\frac{1}{\sqrt{nb_n}}&=&o\left(\left(\frac{1}{\rho_n}\right)^p\right).\label{rate2}%
\end{eqnarray}
The right hand sides of (\ref{rate1}) and (\ref{rate2}) are
increasing and decreasing functions, respectively, of $\rho_n$.
Assumption 3, together with (\ref{2nd}), indicate that $\rho_n$ goes
to infinity as $n$ goes to infinity. The rate given by (\ref{rate1})
will be unduly slow if $\rho_n$ goes to infinity too slowly, while
the rate given by (\ref{rate2}) will be unduly slow if $\rho_n$ goes
to infinity too fast. At either event, $1/\sqrt{nb_n}$ will have a
sub-optimal rate of convergence to zero. It follows that
$1/\sqrt{nb_n}$ has the fastest convergence to zero if
$$O\left(\left(\frac{n}{\rho_n}\right)^{\frac q{1+2q}}\right)=O\left(\rho_n^p\right).$$
This condition requires that $\rho_n=O\left(n^{\frac
q{p+q+2pq}}\right)$. For this rate of $\rho_n$, (\ref{rate1})
implies that
$$b_n=o\left(n^{-\frac {p+q}{p+q+2pq}}\right)\quad \mbox{and }\frac1{\sqrt{nb_n}}=o\left(n^{-\frac
{pq}{p+q+2pq}}\right).$$ This completes the proof.\hfill$\Box$

\section*{References}
\begin{description}
\item {\sc Brillinger, D.~R.} (2001).
{\it Time Series Data Analysis and Theory}. Philadelphia: SIAM.

\item {\sc Brockwell, P.~J.} and {\sc Davis, R.~A.} (1991).
{\it Time Series: Theory and Methods}. New York: Springer-Verlag.

\item {\sc Chen, H., Simpson, D.~G.} and {\sc Ying, Z.} (2000).
Infill asymptotics for a stochastic process model with measurement
error. {\it Statist. Sinica} 10, 141�-156.

\item {\sc Constantine, A.~G.} and {\sc Hall, P.} (1994).
Characterizing surface smoothness via estimation of effective
fractal dimension. {\it J. Roy. Statist. Soc. Ser. B} 56, 96�-113.

\item {\sc Fuentes, M.} (2002).
Spectral methods for nonstationary spatial processes. {\it
Biometrika} 89, 197-�210.

\item {\sc Hall, P., Fisher, N.~I.} and {\sc Hoffmann, B.} (1994).
On the nonparametric estimation of covariance functions. {\it Ann.
Statist.} 22, 2115�-2134.

\item {\sc Hoel, P.~G., Port, S.~C.} and {\sc Stone, C.~J.} (1972).
{\it Introduction to Stochastic Processes}. Boston: Houghton
Mifflin.

\item {\sc Kay, S.~M.} (1999).
{\it Modern Spectral Estimation: Theory and Application}. Englewood
Cliffs, New Jersey: Prentice Hall.

\item {\sc Lahiri, S.~N.} (1999).
Asymptotic distribution of the empirical spatial cumulative
distribution function predictor and prediction bands based on a
subsampling method. {\it Probab. Theory Related Fields} 114, 55�-84.

\item {\sc Leonov, V.~P.} and {\sc Shiryayev, A.~N.} (1959).
On a method of calculation of semi-invariants. {\it Theory Probab.
Appl.} 4, 319�-329.

\item {\sc Lim, C.~Y.} and {\sc Stein, M.} (2008).
Properties of spatial cross-periodograms using fixed-domain
asymptotics. {\it J. Multivariate Anal.} 99, 1962-�1984.

\item {\sc Marvasti, F.~A.} (2001).
{\it Nonuniform Sampling}. New York: Kluwer Plenum.

\item {\sc Masry, E.} (1978).
Alias-free sampling: An alternative conceptualization and its
application. {\it IEEE Trans. Inf. Theor.} IT-24, 173�-183.

\item {\sc Parzen, E.} (1957).
On consistent estimation of the spectrum of the stationary time
series. {\it Annals of Mathematical Statistics} 28, 329�-348.

\item {\sc Resnick, S.~I.} (1987).
{\it Extreme Values, Regular Variation and Point Processes}. New
York : Springer-Verlag.

\item {\sc Shapiro, H.~S.} and {\sc Silverman, R.~A.} (1960).
Alias-free sampling of random noise. {\it J. Soc. Indust. Appl.
Math.} 8, 225-�248.

\item {\sc Shorak, G.~R.} and {\sc Wellner, J.~A.} (1986).
{\it Empirical Processes with Applications to Statistics}. New York:
John Wiley.

\item {\sc Srivastava, R.} and {\sc Sengupta, D.} (2010).
Consistent estimation of non-bandlimited spectral density from
uniformly spaced samples. to appear in {\it IEEE Trans. Inf.
Theor.}, (preprint available in {\small \verb8arXiv:0906.5045v28}).

\item {\sc Stein, M.~L.} (1995). Fixed-domain asymptotics for spatial
periodograms. {\it J. Amer. Statist. Assoc.} 90, 1277�-1288.

\item {\sc Zhang, H.} and {\sc Zimmerman, D.~L.} (2005). Towards
reconciling two asymptotic frameworks in spatial statistics. {\it
Biometrika} 92, 921�-936.
\end{description}

\end{document}